\documentclass[10pt,letterpaper]{article}
\usepackage{amssymb}

\usepackage{graphicx}
\usepackage{amsmath}

\newtheorem{theorem}{Theorem}

\newtheorem{corollary}[theorem]{Corollary}

\newtheorem{definition}[theorem]{Definition}
\newtheorem{example}[theorem]{Example}

\newtheorem{lemma}[theorem]{Lemma}
\newtheorem{notation}[theorem]{Notation}

\newtheorem{proposition}[theorem]{Proposition}
\newtheorem{remark}[theorem]{Remark}

\newenvironment{proof}[1][Proof]{\textbf{#1.} }{\ \rule{0.5em}{0.5em}}

\begin{document}

\title{Morphisms of locally compact groupoids endowed with Haar systems}
\author{M\u{a}d\u{a}lina Roxana Buneci \and Piotr Stachura}
\date{}
\maketitle

\begin{abstract}
We shall generalize the notion of groupoid morphism given by Zakrzewski ( 
\cite{za1}, \cite{za2}) to the setting of locally compact $\sigma $-compact
Hausdorff groupoids endowed with Haar systems. To each groupoid $\Gamma $
endowed with a Haar system $\lambda $ we shall associate a $C^{\ast }$%
-algebra $C^{\ast }\left( \Gamma ,\lambda \right) $, and we construct a
covariant functor $\left( \Gamma ,\lambda \right) \rightarrow C^{\ast
}\left( \Gamma ,\lambda \right) $ from the category of locally compact, $%
\sigma $-compact, Hausdorff groupoids endowed with Haar systems to the
category of $C^{\ast }$-algebras (in the sense of \cite{wo}). If $\Gamma $
is second countable and measurewise amenable, then $C^{\ast }\left( \Gamma
,\lambda \right) $ coincides with the full and the reduced $C^{\ast }$%
-algebras associated to $\Gamma $ and $\lambda .$
\end{abstract}

{\small AMS 2000 Subject Classification: 22A22, 43A22, 46L99.}

{\small Key Words: locally compact groupoid, morphism, $C^{\ast }$-algebra.} 

\section{Introduction}

The purpose of this paper is extend the notion of morphism of groupoids
introduced in \cite{za1,za2} to locally compact $\sigma$-compact groupoids
endowed with Haar systems and to use the extension to construct a covariant
functor from this category to the category of $C^{\ast }$-algebras.

If $\Gamma $ and $G$ are two locally compact, $\sigma $-compact, Hausdorff
groupoids and if $\lambda =\left\{ \lambda ^{u},u\in \Gamma ^{\left(
0\right) }\right\} $ (respectively, $\nu =\left\{ \nu ^{t},t\in G^{\left(
0\right) }\right\} $) is a Haar system on $\Gamma $ (respectively, on $G$),
then by a morphism from $\left( \Gamma ,\lambda \right) $ to $\left( G,\nu
\right) $ we shall mean a left action of $\Gamma $ on $G$, which commutes
with multiplication on $G$, and which satisfies a ''compatibility''
condition with respect to the Haar systems on $\Gamma $ and $G$. This notion
of morphism reduces to a group homomorphism if $\Gamma $ and $G$ are groups,
and to a map in the reverse direction if $\Gamma $ and $G$ are sets.

To each groupoid $\Gamma $ endowed with a Haar system $\lambda $ we shall
associate a $C^{\ast }$-algebra $C^{\ast }\left( \Gamma ,\lambda \right) $
in the following way. We shall consider the space $C_{c}\left( \Gamma
\right) $ of complex-valuated continuous functions with compact support on $%
\Gamma $, which is a topological $\ast $-algebra under the usual convolution
and involution. To each morphism $h$ from $\left( \Gamma ,\lambda \right) $
to a groupoid $\left( G,\nu \right) $ and each unit $t\in G^{\left( 0\right)
}$ we shall associate a $\ast $-representation $\pi _{h,t}$ of $C_{c}\left(
\Gamma \right) $. For any morphism $h$ from $\left( \Gamma ,\lambda \right) $
to a groupoid $\left( G,\nu \right) $ we shall define 
\begin{equation*}
\left\| f\right\| _{h}=\sup\limits_{t}\left\| \pi _{h,t}\left( f\right)
\right\|
\end{equation*}
for all $f\in $ $C_{c}\left( \Gamma \right) $. The $C^{\ast }$-algebra $%
C^{\ast }\left( \Gamma ,\lambda \right) $ will be defined to be the
completion of $C_{c}\left( \Gamma ,\lambda \right) $ in the norm $\left\|
\cdot \right\| =sup_{h}\left\| \cdot \right\| _{h}$, where $h$ runs over all
morphism defined on $\left( \Gamma ,\lambda \right) $. We shall show the
following inequalities 
\begin{equation*}
\left\| f\right\| _{red}\leq \left\| f\right\| \leq \left\| f\right\| _{full}
\end{equation*}
for all $f\in C_{c}\left( \Gamma ,\lambda \right) $, where $\left\| \cdot
\right\| _{red}$ and $\left\| \cdot \right\| _{full}$ are the usual reduced
and full $C^{\ast }$-norms on $C_{c}\left( \Gamma ,\lambda \right) $.
Therefore, according to prop. 6.1.8/p.146 \cite{ar}, if $\left( \Gamma
,\lambda \right) $ is measurewise, amenable then $C^{\ast }\left( \Gamma
,\lambda \right) =C_{full}^{\ast }\left( \Gamma ,\lambda \right)
=C_{red}^{\ast }\left( \Gamma ,\lambda \right) $, where $C_{full}^{\ast
}\left( \Gamma ,\lambda \right) $(respectively, $C_{red}^{\ast }\left(
\Gamma ,\lambda \right) $) is the full (respectively, the reduced) $C^{\ast
} $-algebra associated to $\Gamma $ and $\lambda $. If the principal
associated groupoid of $\Gamma $ is a proper groupoid (but $\Gamma $ is not
necessarily measurewise amenable), then for any quasi invariant measure $\mu 
$ on $\Gamma ^{\left( 0\right) }$ and any $f\in C_{c}\left( \Gamma ,\lambda
\right) $ we shall prove that 
\begin{equation*}
\left\| II_{\mu }\left( f\right) \right\| \leq \left\| f\right\| \text{,}
\end{equation*}
where $II_{\mu }$ is the trivial representation on $\mu $ of $C_{c}\left(
\Gamma ,\lambda \right) $.

Now we are going to establish notation. Relevant definitions can be found in
several places (e.g. \cite{re1}, \cite{mu}). For a groupoid $\Gamma$ the set
of \emph{composable pairs} will be denoted by $\Gamma^{(2)}$ and the \emph{%
multiplication map} by $\Gamma^{(2)}\ni(x,y)\mapsto x y \in \Gamma$; the 
\emph{inverse map} by $\Gamma\ni x \mapsto x^{-1}\in \Gamma$. The set of 
\emph{units} by $\Gamma^{(0)}$ and the \emph{domain} and \emph{range} maps
by $d: \Gamma\ni x \mapsto d(x):=x^{-1}x\in \Gamma^{(0)}$, and $r: \Gamma\ni
x \mapsto d(x):=x x^{-1}\in \Gamma^{(0)}$ respectively.

The fibres of the range and the source maps are denoted $\Gamma
^{u}=r^{-1}\left( \left\{ u\right\} \right) $ and $\Gamma _{v}=d^{-1}\left(
\left\{ v\right\} \right) $, respectively. More generally, given the subsets 
$A$, $B$ $\subset \Gamma ^{\left( 0\right) }$, we define $%
\Gamma^{A}=r^{-1}\left( A\right) $, $\Gamma _{B}=d^{-1}\left( B\right) $ and 
$\Gamma _{B}^{A}=r^{-1}\left( A\right) \cap d^{-1}\left( B\right). $ The 
\textit{reduction} of $\Gamma $ to $A\subset \Gamma ^{\left( 0\right) }$ is $%
\Gamma |_{A}=\Gamma _{A}^{A}$. The relation $u \sim v$ iff $%
\Gamma_{v}^{u}\neq \phi $ is an equivalence relation on $\Gamma ^{\left(
0\right) } $. Its equivalence classes are called \textit{orbits}, the orbit
of a unit $u $ is denoted $\left[ u\right] $ and the quotient space for this
equivalence relation is called the \textit{orbit space} of $G$ and denoted $%
G^{\left( 0\right) }/G$. A groupoid is called \textit{transitive} iff it has
a single orbit, or equivalently, if the map 
\begin{equation*}
\theta :\Gamma \rightarrow \left\{ \left( r\left( x\right) ,d\left( x\right)
\right) ,\,x\in \Gamma \right\} ,\,\,\theta \left( x\right) =\left( r\left(
x\right) ,d\left( x\right) \right)
\end{equation*}
is surjective. A groupoid $\Gamma $ is called \textit{principal}, if the
above map $\theta $ is injective. A subset of $G^{\left( 0\right) }$ is said 
\textit{saturated} if it contains the orbits of its elements. For any subset 
$A$ of $\Gamma ^{\left( 0\right) }$, we denote by $\left[ A\right] $ the
union of the orbits $\left[ u\right] $ for all $u\in A$.

A \textit{topological groupoid} consists of a groupoid $\Gamma $ and a
topology compatible with the groupoid structure. This means that:

$\left( 1\right) \;x\rightarrow x^{-1}\;\left[ :\Gamma \rightarrow \Gamma %
\right] $ is continuous.

$\left( 2\right) $\ $\left( x,y\right) \;\left[ :\Gamma ^{\left( 2\right)
}\rightarrow \Gamma \right] $ is continuous where $\Gamma ^{\left( 2\right)
} $ has the induced topology from $\Gamma \times \Gamma $.

We are exclusively concerned with topological groupoids which are locally
compact Hausdorff.

In the following we will use abbreviation \emph{lcH}-groupoids.

If $\Gamma $ is Hausdorff, then $\Gamma ^{\left( 0\right) }$ is closed in $%
\Gamma $, and $\Gamma ^{\left( 2\right) }$ closed in $\Gamma \times \Gamma $%
.\ It was shown in \cite{ra1} that measured groupoids (in the sense of def.
2.3./p.6 \cite{ha1} ) may be assume to have locally compact topologies, with
no loss of generality.

A subset $A$ of a locally compact groupoid $\Gamma $ is called $r$-\textit{%
(relatively) compact} iff $A\cap r^{-1}\left( K\right) $ is (relatively)
compact for each compact subset $K$ of $\Gamma ^{\left( 0\right) }$.
Similarly, one may define $d$-\textit{(relatively) compact} subsets of $%
\Gamma $. A subset of \ $\Gamma $ which is $r$-(relatively) compact and $d$%
-(relatively) compact is said \textit{conditionally-(relatively) compact}.
If the unit space $\Gamma ^{\left( 0\right) }$ is paracompact, then there
exists a fundamental system of conditionally-(relatively) compact
neighborhoods of $\Gamma ^{\left( 0\right) }$ (see the proof of prop.
II.1.9/p.56 \cite{re1}).

If $X$ is a locally compact space, $C_{c}\left( X\right) $ denotes the space
of complex-valuated continuous functions with compact support. The Borel
sets of a topological space are taken to be the $\sigma $-algebra generated
by open sets.

If $\Gamma $ is a locally compact groupoid, then for each $u\in \Gamma
^{\left( 0\right) }$, $\Gamma _{u}^{u}=\Gamma |_{\left\{ u\right\} }$ is a
locally compact group. We denote by 
\begin{equation*}
\Gamma ^{\prime }=\left\{ x\in \Gamma :\,r\left( x\right) =d\left( x\right)
\right\} =\bigcup_{u\in \Gamma ^{\left( 0\right) }}\Gamma _{u}^{u}
\end{equation*}
the \textit{isotropy group bundle} of $\Gamma $. It is closed in $\Gamma $.

Recall that a (continuous) \textit{Haar system} on a lcH-groupoid $\Gamma$
is a family of positive Radon measures on $\Gamma $, $\lambda =\left\{
\lambda ^{u},\,u\in \Gamma ^{\left( 0\right) }\right\} $, such that

\begin{enumerate}
\item  For all $u\in \Gamma ^{\left( 0\right) }$, $supp(\lambda
^{u})=\,\Gamma ^{u}$.

\item  For all $f\in C_{c}\left( \Gamma \right) $, 
\begin{equation*}
u\rightarrow \int f\left( x\right) d\lambda ^{u}\left( x\right) \,\;\left[
:\Gamma ^{\left( 0\right) }\rightarrow \mathbf{C}\right]
\end{equation*}

is continuous

\item  For all $f\in C_{c}\left( \Gamma \right) $ and all $x\in \Gamma $, 
\begin{equation*}
\int f\left( y\right) d\lambda ^{r\left( x\right) }\left( y\right) \,=\,\int
f\left( xy\right) d\lambda ^{d\left( x\right) }\left( y\right) \text{.}\;
\end{equation*}
\end{enumerate}

Unlike the case of locally compact group, Haar system on groupoid need not
exists, and if it does, it will not usually be unique. The continuity
assumption $2)$ has topological consequences for $\Gamma $. It entails that
the range map $r:\Gamma \rightarrow \Gamma ^{\left( 0\right) }$, and hence
the domain map $d:\Gamma \rightarrow \Gamma ^{\left( 0\right) }$ is open
(prop. I.4 \cite{we}). \emph{Therefore, in this paper we shall always assume
that $r:\Gamma \rightarrow \Gamma ^{\left( 0\right) }$\ is an open map.} For
each $\lambda ^{u}$, \ we denote by $\lambda _{u}$ the image of $\lambda
^{u} $ by the inverse map $x\rightarrow x^{-1}$ (i.e. $\int f\left( y\right)
d\lambda _{u}\left( y\right) =\int f\left( y^{-1}\right) d\lambda ^{u}\left(
y\right) $, \ \ $f\in C_{c}\left( \Gamma \right) $) .

If $\mu $ is a Radon measure on $\Gamma ^{\left( 0\right) }$, then the
measure $\lambda ^{\mu }=\int \lambda ^{u}d\mu \left( u\right) $, defined by 
\begin{equation*}
\int f\left( y\right) d\lambda ^{\mu }\left( y\right) =\int \left( \int
f\left( y\right) d\lambda ^{u}\left( y\right) \right) d\mu \left( u\right) 
\text{, \ \ }f\in C_{c}\left( \Gamma \right)
\end{equation*}
is called the \textit{measure on }$\Gamma $\textit{\ induced by }$\mu $. The
image of $\lambda ^{\mu }$ by the inverse map $x\rightarrow x^{-1}$ is
denoted $\left( \lambda ^{\mu }\right) ^{-1}$. $\mu $ is said \textit{%
quasi-invariant} if its induced measure $\lambda ^{\mu }$ is equivalent to
its inverse $\left( \lambda ^{\mu }\right) ^{-1}$. \ A measure belonging to
the class of a quasi-invariant measure is also quasi-invariant. We say that
the class is \textit{invariant}.

If $\mu $ is a quasi-invariant measure on $\Gamma ^{\left( 0\right) }$ and $%
\lambda ^{\mu }$ is the measure induced on $G$, then the Radon-Nikodym
derivative $\Delta =\frac{d\lambda ^{\mu }}{d\left( \lambda ^{\mu }\right)
^{-1}}$ is called the \textit{modular function} of $\mu $. According to cor.
3.14/p.19 \cite{ha1}, there is a $\mu $-conull Borel subset $U_{0}$ of $%
\Gamma ^{\left( 0\right) }$ such that the restriction of $\Delta $ to $%
\Gamma |_{U_{0}}$ is a homomorphism.

A family of positive Borel measures on $\Gamma $, $\lambda =\left\{ \lambda
^{u},\,u\in \Gamma ^{\left( 0\right) }\right\} $ is called a \textit{Borel
Haar system} if in the definition of the (continuous) Haar system we replace
the condition $2$ by

\begin{description}
\item[$2^{\prime }$]  For all $f\geq 0$ Borel on $\Gamma $, the map 
\begin{equation*}
u\rightarrow \int f\left( x\right) d\lambda ^{u}\left( x\right) \,\;\left[
:\Gamma ^{\left( 0\right) }\rightarrow \mathbf{R}\right]
\end{equation*}
\end{description}

is Borel, and there is a nonnegative Borel function $F$ on $\Gamma $ such
that 
\begin{equation*}
\int F\left( x\right) d\lambda ^{u}\left( x\right) =1\text{ for all }u\text{.%
}
\end{equation*}

For a $C^{\ast }$-algebra A let $\mathcal{B}\left( A\right) $ be the algebra
of bounded linear map acting on $A$. We say that $a\in \mathcal{B}\left(
A\right)$ \emph{is adjointable} iff there is an element $b\in \mathcal{B}%
\left( A\right) $ such that $y^{\ast }\left( a\left( x\right) \right)
=\left( b\left( y\right) \right) ^{\ast }x$ for all $x,y\in A$. The set of
adjointable $a\in \mathcal{B}\left( A\right) $ is the \textit{multiplier
algebra} of $A$, and will be denoted by $M\left( A\right) $. A \textit{%
morphism }from a $C^{\ast } $-algebra $A$ to a $C^{\ast }$-algebra $B$, is a 
$\ast $-homomorphism $\phi :A\rightarrow M\left( B\right) $ such that the
set \ $\phi \left( A\right) B$ is dense in $B$. Such morphism extends
uniquely to a $*$-homomorphism $\hat{\phi}$ from $M\left( A\right) $ to $%
M\left( B\right) $ by\ $\hat{\phi}\left( m\right) \left( \phi \left(
a\right) b\right) =\phi \left( m\left( a\right) \right) b$ for $m\in M\left(
A\right) $, $a\in A$ and $b\in B$. If $\phi _{1}$ is a morphism from $A$ to $%
B$, and $\phi _{2}$ is a morphism from $B$ to $C$, the composition is
defined by $\hat{\phi}_{2}\phi _{1}:A\rightarrow M\left( C\right) $. $%
C^{\ast }$ algebras with above defined morphisms form a $C^{\ast }$-category
(see \cite{wo}).


\section{The category of locally compact groupoids endowed with Haar systems}

\subsection{Definition of morphism}

\bigskip

Let $\Gamma $ and $G$ be two $\sigma $-compact, \emph{lcH}-groupoids. Let $%
\lambda =\left\{ \lambda ^{u},u\in \Gamma ^{\left( 0\right) }\right\} $
(respectively, $\nu =\left\{ \nu ^{t},t\in G^{\left( 0\right) }\right\} $)
be a Haar system on $\Gamma $ (respectively, on $G$).

By a morphism from $\left( \Gamma ,\lambda \right) $ to $\left( G,\nu
\right) $ we mean a left action of $\Gamma $ on $G$, which commutes with
multiplication on $G$, and which satisfies a ''compatibility'' condition
with respect to the Haar systems on $\Gamma $ and $G$.

\begin{definition}
\label{action}Let $\Gamma $ be a groupoid and $X$ be a set. We say $\Gamma $
acts (to the left) on $X$ if there is a map $\ \rho :X\rightarrow \Gamma
^{\left( 0\right) }$ ( called a momentum map) and a map $\left( \gamma
,x\right) \rightarrow \gamma \cdot x$ from 
\begin{equation*}
\Gamma \ast _{\rho }X=\left\{ \left( \gamma ,x\right) :d\left( \gamma
\right) =\rho \left( x\right) \,\right\}
\end{equation*}
to $X$, called (left) action, such that:

\begin{enumerate}
\item  $\rho \left( \gamma \cdot x\right) =r\left( \gamma \right) $ for all $%
\left( \gamma ,x\right) \in \Gamma \ast _{\rho }X$.

\item  $\rho \left( x\right) \cdot x=x$ for all $x\in X$.

\item  If $\left( \gamma _{2},\gamma _{1}\right) \in \Gamma ^{\left(
2\right) }$ and $\left( \gamma _{1},x\right) \in \Gamma \ast _{\rho }X$,
then $\left( \gamma _{2}\gamma _{1}\right) \cdot x=\gamma _{2}\cdot \left(
\gamma _{1}\cdot x\right) $.
\end{enumerate}

If $\Gamma $ is a topological groupoid and $X$ is a topological space, then
we say that a left action is continuous if the mappings $\rho $ and $\left(
\gamma ,x\right) \rightarrow \gamma \cdot x$ are continuous, where $\Gamma
\ast _{\rho }X$ is endowed with the relative product topology coming from $%
\Gamma \times X$.
\end{definition}

The difference with the definition of action in \cite{mu} (def. 2.12/p. 32
and rem. 2.30/p.45) or in \cite{mrw} is that we do not assume that the
momentum map is surjective and open.

The action is called \emph{free} if $\left( \gamma ,x\right) \in \Gamma \ast
_{\rho }X$ and $\gamma \cdot x=x$ implies $\gamma \in \Gamma ^{\left(
0\right) }$.

The continuous action is called \emph{proper} if the map $\left( \gamma
,x\right) \rightarrow \left( \gamma \cdot x,x\right) $ from $\Gamma \ast
_{\rho }X$ to $X\times X$ is proper (i.e. the inverse image of each compact
subset of $X\times X$ is a compact subset of $\Gamma \ast _{\rho }X$).

In the same manner, we define a \emph{right action} of $\Gamma $ on $X$,
using a continuous map $\sigma :X\rightarrow \Gamma ^{\left( 0\right) }$ and
a map $\left( x,\gamma \right) \rightarrow x\cdot \gamma $ from 
\begin{equation*}
X\ast _{\sigma }\Gamma =\left\{ \left( x,\gamma \right) :\sigma \left(
x\right) =r\left( \gamma \right) \,\right\}
\end{equation*}
to $X$.

The simplest example of proper and free action is the case when the groupoid 
$\Gamma $ acts upon itself by either right or left translation
(multiplication).

\begin{definition}
Let $\Gamma _{1},\Gamma _{2}$ be two groupoids and $X$ be set. Let us assume
that $\Gamma _{1}$ acts to the left on $X$ with momentum map $\rho
:X\rightarrow \Gamma _{1}^{\left( 0\right) }$, and that $\Gamma _{2}$ acts
to the right on $X$ with momentum map $\sigma :X\rightarrow \Gamma
_{2}^{\left( 0\right) }$. We say that the action commute if

\begin{enumerate}
\item  $\rho \left( x\cdot \gamma _{2}\right) =$ $\rho \left( x\right) $ for
all $\left( x,\gamma _{2}\right) \in X\ast _{\sigma }\Gamma _{2}$ and $%
\sigma \left( \gamma _{1}\cdot x\right) =\sigma \left( x\right) $ for all $%
\left( \gamma _{1},x\right) \in \Gamma _{1}\ast _{\rho }X.$

\item  $\gamma _{1}\cdot \left( x\cdot \gamma _{2}\right) =\left( \gamma
_{1}\cdot x\right) \cdot \gamma _{2}$ for all $\left( \gamma _{1},x\right)
\in \Gamma _{1}\ast _{\rho }X$, $\left( x,\gamma _{2}\right) \in X\ast
_{\sigma }\Gamma _{2}$.
\end{enumerate}
\end{definition}

\bigskip

\begin{definition}
\label{morph_a} Let $\Gamma $ and $G$ be two groupoids. By an algebraic
morphism from $\Gamma $ to $G$ we mean a left action of $\Gamma $ on $G$
which commutes with the multiplication on $G$.

The morphism is said continuous if the action of $\Gamma $ on $G$ is
continuous (assuming that $\Gamma $ and $G$ are topological spaces).
\end{definition}

Let us note that if we have a morphism in the sense of the preceding
definition and if $\rho :G\rightarrow \Gamma $ is the momentum map of the
left action, then $\rho =\rho \circ r$. Indeed, for any $x\in G$, we have $%
\rho \left( x\right) =\rho \left( xx^{-1}\right) =\rho \left( r\left(
x\right) \right) $ because of the fact that left action of $\Gamma $ on $G$
commutes with the multiplication on $G$.

Therefore an algebraic morphism $h:\Gamma \mbox{$\,$%
\rule[0.5ex]{1.1em}{0.2pt}$\triangleright\,$}G$ is given by two maps

\begin{enumerate}
\item  $\rho _{h}:G^{\left( 0\right) }\rightarrow \Gamma ^{\left( 0\right) }$

\item  $\left( \gamma ,x\right) \rightarrow \gamma \cdot _{h}x$ from $\Gamma
\star _{h}G$ to $G$, where 
\begin{equation*}
\Gamma \star _{h}G=\left\{ \left( \gamma ,x\right) \in \Gamma \times
G:d\left( \gamma \right) =\rho _{h}\left( r\left( x\right) \right) \right\}
\end{equation*}
\end{enumerate}

satisfying the following conditions:

\begin{description}
\item[(1)]  $\rho _{h}\left( r\left( \gamma \cdot _{h}x\right) \right)
=r\left( \gamma \right) $ for all $\left( \gamma ,x\right) \in \Gamma \star
_{h}G$.

\item[(2)]  $\rho _{h}\left( r\left( x\right) \right) \cdot _{h}x=x$ for all 
$x\in G$.

\item[(3)]  $\left( \gamma _{1}\gamma _{2}\right) \cdot _{h}x=\gamma
_{1}\cdot _{h}\left( \gamma _{2}\cdot _{h}x\right) $ for all $\left( \gamma
_{1},\gamma _{2}\right) \in \Gamma ^{\left( 2\right) }$ and all $\left(
\gamma _{2},x\right) \in \Gamma \star _{h}G$.

\item[(4)]  $d\left( \gamma \cdot _{h}x\right) =d\left( x\right) $ for all $%
\left( \gamma ,x\right) \in \Gamma \star _{h}G$.

\item[(5)]  $\left( \gamma \cdot _{h}x_{1}\right) x_{2}=\gamma \cdot
_{h}\left( x_{1}x_{2}\right) $ for all $\left( \gamma ,x_{1}\right) \in
\Gamma \star _{h}G$ and $\left( x_{1},x_{2}\right) \in G^{\left( 2\right) }$.
\end{description}

In the case continuous morphism the map $\rho _{h}$ is a continuous map. The
map $\rho _{h}$ is not necessarily open or surjective. However, the image of 
$\rho _{h}$ is always a saturated subset of $\Gamma ^{\left( 0\right) }$.
Indeed, let $v\sim u=\rho _{h}\left( t\right) $ and let $\gamma \in \Gamma $
be such that $r\left( \gamma \right) =v$ and $d\left( \gamma \right) =u$.
Then $v$ belongs to the image of $\rho _{h}$ because $v=r\left( \gamma
\right) =\rho _{h}\left( r\left( \gamma \cdot _{h}t\right) \right) $.

\begin{remark}
Let $h$ be an algebraic morphism from $\Gamma $ to $G$ ( in the sense of
Def. \ref{morph_a}). Then $h$ is determined by $\rho _{h}$ and the
restriction of the action to 
\begin{equation*}
\left\{ \left( \gamma ,t\right) \in \Gamma \times G^{\left( 0\right)
}:d\left( \gamma \right) =\rho _{h}\left( t\right) \right\} \,\text{.}
\end{equation*}
Indeed, using the condition 5, one obtains 
\begin{equation*}
\gamma \cdot _{h}x=\left( \gamma \cdot _{h}r\left( x\right) \right) x
\end{equation*}
Let us also note that 
\begin{equation*}
\begin{array}{l}
\left( \gamma _{1}\gamma _{2}\right) \cdot _{h}x=\left( \left( \gamma
_{1}\gamma _{2}\right) \cdot _{h}r\left( x\right) \right) x=\gamma _{1}\cdot
_{h}\left( \gamma _{2}\cdot _{h}r\left( x\right) \right) x \\[3mm] 
=\left( \gamma _{1}\cdot _{h}r\left( \gamma _{2}\cdot _{h}r\left( x\right)
\right) \right) \left( \gamma _{2}\cdot _{h}r\left( x\right) \right) x.
\end{array}
\end{equation*}
Consequently, for any $\gamma \in \Gamma $ and any $t\in G^{\left( 0\right)
} $ with $\rho _{h}\left( t\right) =d\left( \gamma \right) $, we have 
\begin{equation*}
\left( \gamma ^{-1}\cdot _{h}r\left( \gamma \cdot _{h}t\right) \right)
\left( \gamma \cdot _{h}t\right) =\left( \gamma ^{-1}\gamma \right) \cdot
_{h}t=d\left( \gamma \right) \cdot _{h}t=\rho _{h}\left( t\right) \cdot
_{h}t=t.
\end{equation*}
Thus for any $\gamma \in \Gamma $ and any $t\in G^{\left( 0\right) }$ with $%
\rho _{h}\left( t\right) =d\left( \gamma \right) $, 
\begin{equation*}
\left( \gamma \cdot _{h}t\right) ^{-1}=\gamma ^{-1}\cdot _{h}r\left( \gamma
\cdot _{h}t\right) .
\end{equation*}
Therefore, algebraically, the notion of morphism in the sense of def. \ref
{morph_a} is the same with that introduced in \cite{za1} (p. 351). In order
to prove the equivalence of these definitions, we can use Prop. 2.7/p. 5\cite
{st}, taking $f=\rho _{h}$ and $g\left( \gamma ,t\right) =\gamma \cdot _{h}t$%
.
\end{remark}

\begin{remark}
Let $h:\Gamma \mbox{$\,$\rule[0.5ex]{1.1em}{0.2pt}$\triangleright\,$}G$ be a
continuous morphism of \emph{lcH}-groupoids ( in the sense of Def. \ref
{morph_a}). Then $G$ is left $\Gamma $-space under the action $\left( \gamma
,x\right) \rightarrow \gamma \cdot _{h}x$, and a right $G$-space under the
multiplication on $G$. $G$ is a \ correspondence in the sense of Def. 2/p.
234\cite{sta} if and only if the left action of $\Gamma $ on $G$ is proper
and $\rho _{h}$ is open and injective. $G$ is a regular bibundle in the
sense of Def. 6/p.103 \cite{la} if and only if the action of $\Gamma $ is
free and transitive along the fibres of $d$ (this means that for all $u\in
G^{\left( 0\right) }$ and $x$ satisfying $d\left( x\right) =x$, there is $%
\gamma \in \Gamma $ such that $\gamma \cdot _{h}u=x$). Therefore, the notion
of morphism introduced in Def. \ref{morph_a} is not cover by the notions
used in \cite{sta} and \cite{la}.
\end{remark}

\begin{lemma}
Let $\Gamma $ and $G$ be two groupoids and let $h$ be an algebraic morphism
from $\Gamma $ to $G$ ( in the sense of Definition \ref{morph_a}). Then the
function 
\begin{equation*}
\left( \gamma ,t\right) \rightarrow \gamma \cdot _{h_{0}}t:=r\left( \gamma
\cdot _{h}t\right)
\end{equation*}
from $\left\{ \left( \gamma ,t\right) \in \Gamma \times G^{\left( 0\right)
}:d\left( \gamma \right) =\rho _{h}\left( t\right) \right\} \,$to $G^{\left(
0\right) }$ defines an action of $\Gamma $ to $G^{\left( 0\right) }$ with
the momentum map $\rho _{h}$.
\end{lemma}

\begin{proof}
Let $\left( \gamma _{1},\gamma _{2}\right) \in \Gamma ^{\left( 2\right) }$
and $\left( \gamma _{1},x\right) \in \left\{ \left( \gamma ,t\right) \in
\Gamma \times G^{\left( 0\right) }:d\left( \gamma \right) =\rho _{h}\left(
t\right) \right\} $. Using the computation in the preceding remark, we
obtain 
\begin{equation*}
\begin{array}{l}
\left( \gamma _{1}\gamma _{2}\right) \cdot _{h_{0}}x=r\left( \left( \gamma
_{1}\gamma _{2}\right) \cdot _{h}x\right) =r\left( \left( \left( \gamma
_{1}\gamma _{2}\right) \cdot _{h}r\left( x\right) \right) x\right) \\[3mm] 
=r\left( \left( \left( \gamma _{1}\gamma _{2}\right) \cdot _{h}r\left(
x\right) \right) \right) =r\left( \gamma _{1}\cdot _{h}\left( \gamma
_{2}\cdot _{h}r\left( x\right) \right) \right) \\ 
=\gamma _{1}\cdot _{h_{0}}\left( \gamma _{2}\cdot _{h_{0}}r\left( x\right)
\right) \text{.}
\end{array}
\end{equation*}
\end{proof}

\begin{notation}
\label{semi}Let $\Gamma $ and $G$ be two groupoids and let $h$ be an
algebraic morphism from $\Gamma $ to $G$ ( in the sense of def. \ref{morph_a}%
). Let us denote 
\begin{eqnarray*}
G\rtimes _{h}\Gamma &=&\left\{ \left( x,\gamma \right) \in G\times \Gamma
:\rho _{h}\left( r\left( x\right) \right) =r\left( \gamma \right) \right\} \\
G^{\left( 0\right) }\rtimes _{h_{0}}\Gamma &=&\left\{ \left( t,\gamma
\right) \in G^{\left( 0\right) }\times \Gamma :\rho _{h}\left( t\right)
=r\left( \gamma \right) \right\}
\end{eqnarray*}
$G\rtimes _{h}\Gamma $, respectively $G^{\left( 0\right) }\rtimes
_{h_{0}}\Gamma $, can be viewed as groupoid under the operations 
\begin{equation*}
\begin{array}{l}
\left( x,\gamma \right) ^{-1}=\left( \gamma ^{-1}\cdot _{h}x,\gamma
^{-1}\right) \\[3mm] 
\left( x,\gamma _{1}\right) \left( \gamma _{1}^{-1}\cdot _{h}x,\gamma
_{2}\right) =\left( x,\gamma _{1}\gamma _{2}\right)
\end{array}
\end{equation*}
respectively, 
\begin{equation*}
\begin{array}{l}
\left( t,\gamma \right) ^{-1}=\left( \gamma ^{-1}\cdot _{h_{0}}t,\gamma
^{-1}\right) \\[3mm] 
\left( t,\gamma _{1}\right) \left( \gamma _{1}^{-1}\cdot _{h_{0}}t,\gamma
_{2}\right) =\left( t,\gamma _{1}\gamma _{2}\right) .
\end{array}
\end{equation*}
(where $\gamma \cdot _{h_{0}}t:=r\left( \gamma \cdot _{h}t\right) $ as in
preceding lemma).
\end{notation}

If $\Gamma $ and $G$ are lcH-groupoids, then $G\rtimes _{h}\Gamma $ and $%
G^{\left( 0\right) }\rtimes _{h_{0}}\Gamma $ are lcH-groupoids. If $\lambda
=\left\{ \lambda ^{u},u\in \Gamma ^{\left( 0\right) }\right\} $ is a Haar
system on $\Gamma $ and if the morphism $h$ is continuous, then $\left\{
\varepsilon _{x}\times \lambda ^{\rho _{h}\left( r\left( x\right) \right)
},x\in G\right\} $ is a Haar system on $G\rtimes _{h}\Gamma $ (where $%
\varepsilon _{x}$ is the unit point mass at $x$) and $\left\{ \varepsilon
_{t}\times \lambda ^{\rho _{h}\left( t\right) },t\in G^{\left( 0\right)
}\right\} $ is a Haar system on $G^{\left( 0\right) }\rtimes _{h_{0}}\Gamma $%
.

\begin{definition}
\label{morph}Let $\Gamma $ and $G$ be two $\sigma $-compact, \emph{lcH}%
-groupoids. Let $\lambda =\left\{ \lambda ^{u},u\in \Gamma ^{\left( 0\right)
}\right\} $ (respectively, $\nu =\left\{ \nu ^{t},t\in G^{\left( 0\right)
}\right\} $) be a Haar system on $\Gamma $ (respectively, on $G$). By a
morphism $h:\left( \Gamma ,\lambda \right) \mbox{$\,$%
\rule[0.5ex]{1.1em}{0.2pt}$\triangleright\,$}\left( G,\nu \right) $ we mean
a continuous morphism from $\Gamma $ to $G$ (in the sense of def.\ref
{morph_a}) which satisfies in addition the following condition:

\begin{description}
\item[(6)]  \label{quasi} There exists a continuous positive function $%
\Delta _{h}$ on 
\begin{equation*}
G\rtimes _{h}\Gamma =\left\{ \left( x,\gamma \right) \in G\times \Gamma
:\rho _{h}\left( r\left( x\right) \right) =r\left( \gamma \right) \right\}
\end{equation*}
such that 
\begin{equation*}
\int \int g(\gamma ^{-1}\cdot _{h}x,\gamma ^{-1})\Delta _{h}(\gamma
^{-1}\cdot _{h}x,\gamma ^{-1})d\lambda ^{\rho _{h}\left( r\left( x\right)
\right) }\left( \gamma \right) d\nu _{t}\left( x\right)
\end{equation*}
\begin{equation*}
=\int \int g(\gamma ,x)d\lambda ^{\rho _{h}\left( r\left( x\right) \right)
}\left( \gamma \right) d\nu _{t}\left( x\right)
\end{equation*}
for all $t\in G^{\left( 0\right) }$ and all Borel nonnegative functions $g$
on $G\rtimes _{h}\Gamma $.
\end{description}
\end{definition}

\begin{remark}
The condition \ref{quasi} in the preceding definition means that each
measure $\nu _{t}$ is quasi-invariant with respect to the Haar system $%
\left\{ \varepsilon _{x}\times \lambda ^{\rho _{h}\left( r\left( x\right)
\right) },x\in G\right\} $ on $G\rtimes _{h}\Gamma $.
\end{remark}

\begin{example}
\label{id}Let $\Gamma $ be a $\sigma $-compact, lcH-groupoid, endowed with a
Haar system $\lambda =\left\{ \lambda ^{u},u\in \Gamma ^{\left( 0\right)
}\right\} $. Let us define a morphism $l:\left( \Gamma ,\lambda \right) %
\mbox{$\,$\rule[0.5ex]{1.1em}{0.2pt}$\triangleright\,$} \left( \Gamma
,\lambda \right) $ by $\rho _{l}=id_{\Gamma ^{\left( 0\right) }}$ and $%
\gamma \cdot _{l}x=\gamma x$ (multiplication on $\Gamma $). It is easy to
check that the conditions in the Definition \ref{morph} are satisfied with $%
\Delta _{l}\equiv 1$.
\end{example}

\begin{lemma}
Let $h:\left( \Gamma ,\lambda \right) \mbox{$\,$\rule[0.5ex]{1.1em}{0.2pt}$%
\triangleright\,$}\left( G,\nu \right) $ be a morphism of $\sigma $-compact,
lcH-groupoids with Haar systems. Then the function $\Delta _{h}$ that
appears in the condition \ref{quasi} of the Definition \ref{morph} satisfies 
$\Delta _{h}\left( x,\gamma \right) =$ $\Delta _{h}\left( r\left( x\right)
,\gamma \right) $ for all $\left( x,\gamma \right) \in G\rtimes _{h}\Gamma $.
\end{lemma}

\begin{proof}
Let $f\geq 0$ be a Borel function on $G\rtimes _{h}\Gamma $. For each $t\in
G^{\left( 0\right) }$ and each $x_{0}\in G^{t}$, we have 
\begin{eqnarray*}
&&\int f\left( x,\gamma \right) d\lambda ^{\rho _{h}\left( r\left( x\right)
\right) }\left( \gamma \right) d\nu _{t}\left( x\right) \\
&=&\int f\left( xx_{0}^{-1},\gamma \right) d\lambda ^{\rho _{h}\left(
r\left( x\right) \right) }\left( \gamma \right) d\nu _{d\left( x_{0}\right)
}\left( x\right) \\
&=&\int f\left( \gamma ^{-1}\cdot _{h}xx_{0}^{-1},\gamma ^{-1}\right) \Delta
_{h}\left( \gamma ^{-1}\cdot _{h}x,\gamma ^{-1}\right) d\lambda ^{\rho
_{h}\left( r\left( x\right) \right) }\left( \gamma \right) d\nu _{d\left(
x_{0}\right) }\left( x\right) \\
&=&\int f\left( \left( \gamma ^{-1}\cdot _{h}r\left( x\right) \right)
xx_{0}^{-1},\gamma ^{-1}\right) \Delta _{h}\left( \left( \gamma ^{-1}\cdot
_{h}r\left( x\right) \right) x,\gamma ^{-1}\right) d\lambda ^{\rho
_{h}\left( r\left( x\right) \right) }\left( \gamma \right) d\nu _{d\left(
x_{0}\right) }\left( x\right) \\
&=&\int f\left( \left( \gamma ^{-1}\cdot _{h}r\left( x\right) \right)
x,\gamma ^{-1}\right) \Delta _{h}\left( \left( \gamma ^{-1}\cdot _{h}r\left(
x\right) \right) xx_{0},\gamma ^{-1}\right) d\lambda ^{\rho _{h}\left(
r\left( x\right) \right) }\left( \gamma \right) d\nu _{t}\left( x\right) \\
&=&\int f\left( \gamma ^{-1}\cdot _{h}x,\gamma ^{-1}\right) \Delta
_{h}\left( \gamma ^{-1}\cdot _{h}xx_{0},\gamma ^{-1}\right) d\lambda ^{\rho
_{h}\left( r\left( x\right) \right) }\left( \gamma \right) d\nu _{t}\left(
x\right) \\
&=&\int f\left( x,\gamma \right) \Delta _{h}\left( xx_{0},\gamma \right)
\Delta _{h}\left( \gamma ^{-1}\cdot _{h}x,\gamma ^{-1}\right) d\lambda
^{\rho _{h}\left( r\left( x\right) \right) }\left( \gamma \right) d\nu
_{t}\left( x\right)
\end{eqnarray*}

Thus for all $t\in G^{\left( 0\right) }$ and $x_{0}\in G^{t}$and almost all $%
\left( x,\gamma \right) \in G\rtimes _{h}\Gamma $, 
\begin{eqnarray*}
1 &=&\Delta _{h}\left( xx_{0},\gamma \right) \Delta _{h}\left( \gamma
^{-1}\cdot _{h}x,\gamma ^{-1}\right) \\
&=&\Delta _{h}\left( xx_{0},\gamma \right) \Delta _{h}\left( x,\gamma
\right) ^{-1}.
\end{eqnarray*}
Therefore $\Delta _{h}\left( xx_{0},\gamma \right) =\Delta _{h}\left(
x,\gamma \right) $ for $\int \lambda ^{\rho _{h}\left( r\left( x\right)
\right) }\left( \gamma \right) d\nu _{t}\left( x\right) $ -a.a.$\left(
x,\gamma \right) \in G\rtimes _{h}\Gamma $. Since $\Delta _{h}$ is a
continuous function and $\int \lambda ^{\rho _{h}\left( r\left( x\right)
\right) }\left( \gamma \right) d\nu _{t}\left( x\right) $ is a measure of
full support on 
\begin{equation*}
\left\{ \left( x,\gamma \right) \in G\rtimes _{h}\Gamma :\,d\left( x\right)
=t\right\} \text{,}
\end{equation*}
it follows that $\Delta _{h}\left( xx_{0},\gamma \right) =\Delta _{h}\left(
x,\gamma \right) $ for all $\left( x,\gamma \right) \in G\rtimes _{h}\Gamma $
with $d\left( x\right) =r\left( x_{0}\right) $. Particularly, for $%
x_{0}=x^{-1}$, it follows that $\Delta _{h}\left( r\left( x\right) ,\gamma
\right) =\Delta _{h}\left( x,\gamma \right) $.
\end{proof}

We shall prove that for particular classes of groupoids we can choose a Haar
system on $G$ such that the condition \ref{quasi} is satisfied. In order to
do this we need some results on the structure of the Haar systems, as
developed by J. Renault in Section $1$ of \cite{re2} and also by A. Ramsay
and M.E. Walter in Section $2$ of \cite{rawm}. In \cite{re2} Jean Renault
constructs a Borel Haar system for $G^{\prime }$(the isotropy group bundle
of a locally compact groupoid $G$ which has a fundamental system of
conditionally-compact neighborhoods of $G^{\left( 0\right) }$). One way to
do this is to choose a function $F_{0}$ continuous with conditionally
support which is nonnegative and equal to $1$ at each $t\in G^{\left(
0\right) }.$ Then for each $t\in G^{\left( 0\right) }$ choose a left Haar
measure $\beta _{t}^{t}$ on $G_{t}^{t}$ so the integral of $F_{0}$ with
respect to $\beta _{t}^{t}$ is $1$. If the restriction of $r$ to $G^{\prime
} $is open, then $\left\{ \beta _{t}^{t},\,t\in G^{\left( 0\right) }\right\} 
$ is a Haar system for $G^{\prime }$(Lemma 1.3/p. 6 \cite{re2}).

Renault defines $\beta _{s}^{t}=x\beta _{s}^{s}$ if $x\in G_{s}^{t}$ (where $%
x\beta _{s}^{s}\left( f\right) =\int f\left( xy\right) d\beta _{s}^{s}\left(
y\right) $). If $z$ is another element in $G_{s}^{t}$, then $x^{-1}z\in
G_{s}^{s}$, and since $\beta _{s}^{s}$ is a left Haar measure on $G_{s}^{s}$%
, it follows that $\beta _{s}^{t}$ is independent of the choice of $x$. If $K
$ is a compact subset of $G$, then $\sup\limits_{t,s}\beta _{s}^{t}\left(
K\right) <\infty $. Renault also defines a $1$-cocycle $\delta _{G}$ on $G$
such that for every $t\in G^{\left( 0\right) }$, $\delta |_{G_{t}^{t}}$ is
the modular function for $\beta _{t}^{t}$. With this apparatus in place,
Renault describes a decomposition of the Haar system $\left\{ \nu
^{t},\,t\in G^{\left( 0\right) }\right\} $ for $G$ over the equivalence
relation $R$ (the principal groupoid associated to $G$). He proves that
there is a unique Borel Haar system $\tilde{\nu}$ for $R$ \ with the
property that 
\begin{equation*}
\nu ^{t}=\int \beta _{s}^{q}d\tilde{\nu}^{t}\left( q,s\right) \text{ \ for
all }t\in G^{\left( 0\right) }\text{.}
\end{equation*}
In Section $2$ \cite{rawm} A. Ramsay and M.E. Walter prove that 
\begin{equation*}
\sup\limits_{t}\tilde{\nu}^{t}\left( \left( r,d\right) \left( K\right)
\right) <\infty \text{, for all compact }K\subset G
\end{equation*}
For each $t\in G^{\left( 0\right) }$ the measure $\tilde{\nu}^{t}$ is
concentrated on $\left\{ t\right\} \times \left[ t\right] $. Therefore there
is a measure $\tilde{\beta}^{t}$ concentrated on $\left[ t\right] $ such
that $\tilde{\nu}^{t}=\varepsilon _{t}\times \tilde{\beta}^{t}$, where $%
\varepsilon _{t}$ is the unit point mass at $t$. Since $\left\{ \tilde{\nu}%
^{t},t\in G^{\left( 0\right) }\right\} $ is a Haar system, we have $\tilde{%
\beta}^{t}=\tilde{\beta}^{s}$ for all $\left( t,s\right) \in R$, and the
function 
\begin{equation*}
t\rightarrow \int f\left( s\right) \tilde{\beta}^{t}\left( s\right) 
\end{equation*}
is Borel for all $f\geq 0$ Borel on $G^{\left( 0\right) }$. If $\mu $ is a
quasi-invariant measure for $\left\{ \nu ^{t},t\in G^{\left( 0\right)
}\right\} $, then $\mu $ is a quasi-invariant measure for $\left\{ \tilde{\nu%
}^{t},t\in G^{\left( 0\right) }\right\} $. Also if $\Delta _{\mu ,R}$ is the
modular function associated to $\left\{ \tilde{\nu}^{t},t\in G^{\left(
0\right) }\right\} $ and $\mu $, then $\Delta _{\mu }=\delta _{G}\Delta
_{\mu ,R}\circ \left( r,d\right) $ can serve as the modular function
associated to $\left\{ \nu ^{t},t\in G^{\left( 0\right) }\right\} $ and $\mu 
$. For each $t$ the measure $\tilde{\beta}^{t}$ is quasi-invariant (Section $%
2$ \cite{rawm}). It is easy to see that $\Delta _{\tilde{\beta}^{t},R}=1$,
and consequently $\Delta _{\tilde{\beta}^{t}}=\delta _{G}$. Thus if $G$ is a
transitive groupoid $\tilde{\beta}^{t}$ is a quasi-invariant measure of full
support having a continuous modular function ($\Delta _{\tilde{\beta}%
^{t}}=\delta _{G}$). \ More generally, let us assume that the associated
principal groupoid $R$ associated to $G$ is proper. This means that $R$ is a
closed subset of $G^{\left( 0\right) }\times G^{\left( 0\right) }$ endowed
with product topology (or equivalently, $G^{\left( 0\right) }/G$ is a
Hausdorff space) and the map 
\begin{equation*}
\left( r,d\right) :G\rightarrow R,\,\left( r,d\right) \left( x\right)
=\left( r\left( x\right) ,d\left( x\right) \right) 
\end{equation*}
is an open map when $R$ is endowed with the relative product topology coming
from $G^{\left( 0\right) }\times G^{\left( 0\right) }$.  If $\mu $ is a
quasi-invariant Radon measure for the Haar system $\left\{ \nu ^{t},\,t\in
G^{\left( 0\right) }\right\} $, then $\mu _{1}=\int \tilde{\beta}^{t}d\mu
\left( t\right) $ is a Radon measure which is equivalent to $\mu $ (see
Remark 6/p. 232 \cite{bu})$.$ It is easy to prove that $\mu _{1}$ has a
continuous modular function $\Delta _{\mu _{1}}=\delta _{G}$.

We shall call the pair of the system of measures 
\begin{equation*}
\left( \left\{ \beta _{s}^{t}\right\} _{\left( t,s\right) \in R},\left\{ 
\tilde{\beta}^{\dot{t}}\right\} _{\dot{t}\in G^{\left( 0\right) }/G}\right) 
\end{equation*}

(described above) \emph{the decomposition of the Haar system }$\left\{ \nu
^{t},t\in G^{\left( 0\right) }\right\} $\emph{\ over the principal groupoid
associated to }$G$. Also we shall call $\delta _{G}$ the $1$\emph{-cocycle
associated to the decomposition}.

\begin{proposition}
\label{q}Let $\Gamma $ and $G$ be two $\sigma $-compact, lcH-groupoids. Let $%
\lambda =\left\{ \lambda ^{u},u\in \Gamma ^{\left( 0\right) }\right\} $
(respectively, $\nu =\left\{ \nu ^{t},t\in G^{\left( 0\right) }\right\} $)
be a Haar system on $\Gamma $ (respectively, on $G$). Let$\left( \left\{
\beta _{s}^{t}\right\} _{\left( t,s\right) \in R},\left\{ \tilde{\beta}%
^{t}\right\} _{\dot{t}\in G^{\left( 0\right) }/G}\right) $ be the\emph{\ }%
decomposition of the Haar system $\left\{ \nu ^{t},t\in G^{\left( 0\right)
}\right\} $ over the principal groupoid associated to $G$. Let $h$ be
continuous morphism from $\Gamma $ to $G$ (in the sense of Definition \ref
{morph_a}). If there is a continuous positive function $\Delta :G^{\left(
0\right) }\rtimes _{h_{0}}\Gamma \rightarrow \mathbf{R}$, such that $\Delta $
is the modular function of $\tilde{\beta}^{t}$ with respect to the Haar
system $\left\{ \varepsilon _{t}\times \lambda ^{\rho _{h}\left( t\right)
},t\in G^{\left( 0\right) }\right\} $ on $G^{\left( 0\right) }\rtimes
_{h_{0}}\Gamma $ for each $\dot{t}\in G^{\left( 0\right) }/G$, then $%
h:\left( \Gamma ,\lambda \right) \mbox{$\,$\rule[0.5ex]{1.1em}{0.2pt}$%
\triangleright\,$} \left( G,\nu \right) $ is a morphism in the sense of
Definition \ref{morph}).
\end{proposition}

\begin{proof}
\bigskip Let $\delta _{G}$ be the $1$-cocycle associated to the
decomposition $\left( \left\{ \beta _{s}^{t}\right\} _{\left( t,s\right) \in
R},\left\{ \tilde{\beta}^{t}\right\} _{\dot{t}\in G^{\left( 0\right)
}/G}\right) $. Let $g$:$G\rtimes _{h}\Gamma \rightarrow \mathbf{R}$ be Borel
nonnegative function. Then we have 
\begin{equation*}
\begin{array}{l}
\int f\left( \gamma ^{-1}\cdot _{h}x,\gamma ^{-1}\right) d\lambda ^{\rho
_{h}\left( r\left( x\right) \right) }\left( \gamma \right) d\nu _{t}\left(
x\right) \\ 
=\int f\left( \gamma ^{-1}\cdot _{h}x^{-1},\gamma ^{-1}\right) d\lambda
^{\rho _{h}\left( d\left( x\right) \right) }\left( \gamma \right) d\nu
^{t}\left( x\right) \\ 
=\int f\left( \gamma ^{-1}\cdot _{h}x^{-1},\gamma ^{-1}\right) d\lambda
^{\rho _{h}\left( d\left( x\right) \right) }\left( \gamma \right) d\beta
_{s}^{t}\left( x\right) d\tilde{\beta}^{\dot{t}}\left( s\right) \\ 
=\int f\left( \gamma ^{-1}\cdot _{h}x^{-1},\gamma ^{-1}\right) d\beta
_{s}^{t}\left( x\right) d\lambda ^{\rho _{h}\left( s\right) }\left( \gamma
\right) d\tilde{\beta}^{\dot{t}}\left( s\right) \\ 
=\int f\left( \gamma ^{-1}\cdot _{h}x,\gamma ^{-1}\right) \delta _{G}\left(
x\right) ^{-1}d\beta _{t}^{s}\left( x\right) d\lambda ^{\rho _{h}\left(
s\right) }\left( \gamma \right) d\tilde{\beta}^{\dot{t}}\left( s\right) \\ 
=\int f\left( \left( \gamma ^{-1}\cdot _{h}r\left( x\right) \right) x,\gamma
^{-1}\right) \delta _{G}\left( x\right) ^{-1}d\beta _{t}^{s}\left( x\right)
d\lambda ^{\rho _{h}\left( s\right) }\left( \gamma \right) d\tilde{\beta}^{%
\dot{t}}\left( s\right) \\ 
=\int f\left( x,\gamma ^{-1}\right) \delta _{G}\left( \left( \gamma
^{-1}\cdot _{h}s\right) ^{-1}x\right) ^{-1}d\beta _{t}^{r\left( \gamma
^{-1}\cdot _{h}s\right) }\left( x\right) d\lambda ^{\rho _{h}\left( s\right)
}\left( \gamma \right) d\tilde{\beta}^{\dot{t}}\left( s\right) \\ 
=\int f\left( x,\gamma ^{-1}\right) \delta _{G}\left( \left( \gamma \cdot
_{h}r\left( \gamma ^{-1}\cdot _{h}s\right) \right) x\right) ^{-1}d\beta
_{t}^{r\left( \gamma ^{-1}\cdot _{h}s\right) }\left( x\right) d\lambda
^{\rho _{h}\left( s\right) }\left( \gamma \right) d\tilde{\beta}^{\dot{t}%
}\left( s\right) \\ 
=\int f\left( x,\gamma \right) \delta _{G}\left( \left( \gamma ^{-1}\cdot
_{h}s\right) x\right) ^{-1}d\beta _{t}^{s}\left( x\right) \Delta \left(
r\left( \gamma ^{-1}\cdot _{h}s\right) ,\gamma ^{-1}\right) d\lambda ^{\rho
_{h}\left( s\right) }\left( \gamma \right) d\tilde{\beta}^{\dot{t}}\left(
s\right) \\ 
=\int f\left( x^{-1},\gamma \right) \delta _{G}\left( \left( \gamma
^{-1}\cdot _{h}s\right) \right) ^{-1}d\beta _{s}^{t}\left( x\right) \Delta
\left( r\left( \gamma ^{-1}\cdot _{h}s\right) ,\gamma ^{-1}\right) d\lambda
^{\rho _{h}\left( s\right) }\left( \gamma \right) d\tilde{\beta}^{\dot{t}%
}\left( s\right) \\ 
=\int f\left( x^{-1},\gamma \right) \delta _{G}\left( \left( \gamma
^{-1}\cdot _{h}d\left( x\right) \right) \right) ^{-1}\Delta \left( r\left(
\gamma ^{-1}\cdot _{h}d\left( x\right) \right) ,\gamma ^{-1}\right) d\lambda
^{\rho _{h}\left( s\right) }\left( \gamma \right) d\nu ^{t}\left( x\right)
\\ 
=\int f\left( x,\gamma \right) \delta _{G}\left( \left( \gamma ^{-1}\cdot
_{h}r\left( x\right) \right) \right) ^{-1}\Delta \left( r\left( \gamma
^{-1}\cdot _{h}r\left( x\right) \right) ,\gamma ^{-1}\right) d\lambda ^{\rho
_{h}\left( s\right) }\left( \gamma \right) d\nu _{t}\left( x\right) .
\end{array}
\end{equation*}

The condition \ref{quasi} in the Definition \ref{morph} is satisfied taking 
\begin{equation*}
\Delta _{h}\left( x,\gamma \right) =\delta _{G}\left( \gamma ^{-1}\cdot
_{h}r\left( x\right) \right) \Delta \left( r\left( x\right) ,\gamma \right) 
\text{.}
\end{equation*}
\end{proof}

\begin{corollary}
\label{t}Let $\Gamma $ and $G$ be two $\sigma $-compact, \emph{lcH}%
-groupoids. Let $\lambda =\left\{ \lambda ^{u},u\in \Gamma ^{\left( 0\right)
}\right\} $ be a Haar system on $\Gamma $. Let $h$ be continuous morphism
from $\Gamma $ to $G$ (in the sense of Definition \ref{morph_a}). If $G$ is
transitive and there is a quasi-invariant measure with respect to the Haar
system $\left\{ \varepsilon _{t}\times \lambda ^{\rho _{h}\left( t\right)
},t\in G^{\left( 0\right) }\right\} $ on $G^{\left( 0\right) }\rtimes
_{h_{0}}\Gamma $ having the support $G^{\left( 0\right) }$ and continuous
modular function, then we can choose a Haar system $\nu $ on $G$ such that $%
h:\left( \Gamma ,\lambda \right) \mbox{$\,$\rule[0.5ex]{1.1em}{0.2pt}$%
\triangleright\,$}\left( G,\nu \right) $ is a morphism in the sense of
Definition \ref{morph}).
\end{corollary}

\begin{proof}
Let $\tilde{\beta}$ be a quasi-invariant measure with respect to the Haar
system $\left\{ \varepsilon _{t}\times \lambda ^{\rho _{h}\left( t\right)
},t\in G^{\left( 0\right) }\right\} $ on $G^{\left( 0\right) }\rtimes
_{h_{0}}\Gamma $ having $supp(\tilde{\beta})=G^{\left( 0\right) }$ and
continuous modular function. Then 
\begin{equation*}
\left\{ \int \beta _{s}^{t}d\tilde{\beta}\left( s\right) ,\,t\in G^{\left(
0\right) }\right\} \,
\end{equation*}
is Haar system on $G$ satisfying the hypothesis of Proposition \ref{q}.
\end{proof}

\begin{remark}
If the associated principal groupoid of $G^{\left( 0\right) }\rtimes
_{h_{0}}\Gamma $ is proper, then there is a quasi-invariant measure with
respect to the Haar system $\left\{ \varepsilon _{t}\times \lambda ^{\rho
_{h}\left( t\right) },t\in G^{\left( 0\right) }\right\} $, having the
support $G^{\left( 0\right) }$ and continuous modular function.\ The
associated principal groupoid of $G^{\left( 0\right) }\rtimes _{h_{0}}\Gamma 
$ is proper if and only if the set
\begin{equation*}
R_{\Gamma ,G^{\left( 0\right) }}=\left\{ \left( t,r\left( \gamma
^{-1}t\right) \right) :\,t\in G^{\left( 0\right) },\gamma \in \Gamma
,\,r\left( \gamma \right) =\rho _{h}\left( t\right) \right\} 
\end{equation*}
is a closed in $G^{\left( 0\right) }\times G^{\left( 0\right) }$ (endowed
with the product topology) and the map 
\begin{equation*}
\left( t,\gamma \right) \rightarrow \left( t,r\left( \gamma ^{-1}t\right)
\right) 
\end{equation*}
is an open map from $G^{\left( 0\right) }\rtimes _{h_{0}}\Gamma $ to $%
R_{\Gamma ,G^{\left( 0\right) }}$, when $R_{\Gamma ,G^{\left( 0\right) }}$
is endowed with the relative product topology coming from $G^{\left(
0\right) }\times G^{\left( 0\right) }$.
\end{remark}


\subsection{Composition of morphisms}

\begin{definition}
\label{comp} Let $h :\left(\Gamma, \lambda\right)\mbox{$\,$%
\rule[0.5ex]{1.1em}{0.2pt}$\triangleright\,$} \left(G_{1}, \nu\right)$ and $%
k :\left(G_{1}, \lambda\right)\mbox{$\,$\rule[0.5ex]{1.1em}{0.2pt}$%
\triangleright\,$} \left(G_{2}, \eta\right)$ be two morphism of locally
compact groupoids endowed with Haar systems. Let $kh :\left(\Gamma,
\lambda\right)\mbox{$\,$\rule[0.5ex]{1.1em}{0.2pt}$\triangleright\,$}
\left(G_{2}, \eta\right)$ be defined by

\begin{enumerate}
\item  $\rho _{kh}:G_{2}^{\left( 0\right) }\rightarrow \Gamma ^{\left(
0\right) }$ 
\begin{equation*}
\rho _{kh}\left( x_{2}\right) =\rho _{h}\left( \rho _{k}\left( x_{2}\right)
\right) \text{ for all }x_{2}\in G_{2}.
\end{equation*}

\item  $\left( \gamma ,x_{2}\right) \rightarrow \gamma \cdot
_{kh}x_{2}:=\left( \gamma \cdot _{h}\rho _{k}\left( r\left( x_{2}\right)
\right) \right) \cdot _{k}x_{2}$ from $\Gamma \star _{kh}G_{2}$ to $G_{2}$,
where 
\begin{equation*}
\Gamma \star _{kh}G_{2}=\left\{ \left( \gamma ,x_{2}\right) \in \Gamma
\times G_{2}:d\left( \gamma \right) =\rho _{kh}\left( r\left( x_{2}\right)
\right) \right\}
\end{equation*}
\end{enumerate}
\end{definition}

\begin{remark}
Let $h:\left( \Gamma ,\lambda \right) \mbox{$\,$\rule[0.5ex]{1.1em}{0.2pt}$%
\triangleright\,$} \left( G_{1},\nu \right) $ and $k:\left( G_{1},\lambda
\right) \mbox{$\,$\rule[0.5ex]{1.1em}{0.2pt}$\triangleright\,$} \left(
G_{2},\eta \right) $ be two morphism of locally compact groupoids endowed
with Haar systems. Let $kh:\left( \Gamma ,\lambda \right) %
\mbox{$\,$\rule[0.5ex]{1.1em}{0.2pt}$\triangleright\,$} \left( G_{2},\eta
\right) $ be as in Definition \ref{comp}. Then for all $\gamma \in \Gamma $,
all $x_{1}\in G_{1}$ with $\rho _{h}\left( r\left( x_{1}\right) \right)
=d\left( \gamma \right) $ and all $x_{2}\in G_{2}$ with $\rho _{k}\left(
r\left( x_{2}\right) \right) =d\left( x_{1}\right) $, we have 
\begin{equation*}
\begin{array}{l}
\gamma \cdot _{kh}\left( x_{1}\cdot _{k}x_{2}\right) =\left( \gamma \cdot
_{h}\rho _{k}\left( r\left( x_{1}\cdot _{k}x_{2}\right) \right) \right)
\cdot _{k}\left( x_{1}\cdot _{k}x_{2}\right) =\left( \gamma \cdot
_{h}r\left( x_{1}\right) x_{1}\right) \cdot _{k}x_{2} \\[3mm] 
=\left( \gamma \cdot _{h}x_{1}\right) \cdot _{k}x_{2}
\end{array}
\end{equation*}
\end{remark}

\begin{lemma}
Let $(\Gamma,\lambda) $, $(G_{1},\nu) $ and $(G_{2},\eta)$ be $\sigma $%
-compact, lcH-groupoids. If $h:\left( \Gamma ,\lambda \right) %
\mbox{$\,$\rule[0.5ex]{1.1em}{0.2pt}$\triangleright\,$} \left( G_{1},\nu
\right) $ and $k:\left( G_{1},\lambda \right) \mbox{$\,$%
\rule[0.5ex]{1.1em}{0.2pt}$\triangleright\,$} \left( G_{2},\eta \right) $
are morphisms, then $kh:\left( \Gamma ,\lambda \right) \mbox{$\,$%
\rule[0.5ex]{1.1em}{0.2pt}$\triangleright\,$} \left( G_{2},\eta \right) $ is
a morphism.
\end{lemma}

\begin{proof}
Let us check the conditions 3, 5 and 6 in the def. \ref{morph}. For all $%
\left( \gamma _{1},\gamma _{2}\right) \in \Gamma ^{\left( 2\right) }$ and
all $x_{2}\in G_{2}$ with $\left( \gamma _{2},x\right) \in \Gamma \star
_{kh}G_{2}$, we have 
\begin{equation*}
\begin{array}{rl}
\left( \gamma _{1}\gamma _{2}\right) \cdot _{kh}x_{2}= & \left( \left(
\gamma _{1}\gamma _{2}\right) \cdot _{h}\rho _{k}\left( r\left( x_{2}\right)
\right) \right) \cdot _{k}x_{2} \\[3mm] 
= & \left( \gamma _{1}\cdot _{h}r\left( \gamma _{2}\cdot _{h}\rho _{k}\left(
r\left( x_{2}\right) \right) \right) \right) \left( \gamma _{2}\cdot
_{h}\rho _{k}\left( r\left( x_{2}\right) \right) \right) \cdot _{k}x_{2} \\%
[3mm] 
= & \left( \left( \gamma _{1}\cdot _{h}r\left( \gamma _{2}\cdot _{h}\rho
_{k}\left( r\left( x_{2}\right) \right) \right) \right) \right) \cdot
_{k}\left( \left( \gamma _{2}\cdot _{h}\rho _{k}\left( r\left( x_{2}\right)
\right) \right) \cdot _{k}x_{2}\right) \\[3mm] 
= & \left( \left( \gamma _{1}\cdot _{h}r\left( \gamma _{2}\cdot _{h}\rho
_{k}\left( r\left( x_{2}\right) \right) \right) \right) \right) \cdot
_{k}\left( \gamma _{2}\cdot _{kh}x_{2}\right) \\[3mm] 
= & \left( \left( \gamma _{1}\cdot _{h}\rho _{k}\left( r\left( \gamma
_{2}\cdot _{h}\rho _{k}\left( r\left( x_{2}\right) \right) \right) \right)
\right) \cdot _{k}x\right) \cdot _{k}\left( \gamma _{2}\cdot
_{kh}x_{2}\right) \\[3mm] 
= & \gamma _{1}\cdot _{kh}\left( \gamma _{2}\cdot _{kh}x_{2}\right) .
\end{array}
\end{equation*}
For all $\left( \gamma ,x_{2}\right) \in \Gamma \star _{kh}G_{2}$ and $%
\left( x_{2},y_{2}\right) \in G_{2}^{\left( 2\right) }$, we have 
\begin{equation*}
\begin{array}{l}
\left( \gamma \cdot _{kh}x_{2}\right) y_{2}=\left( \left( \gamma \cdot
_{h}\rho _{k}\left( r\left( x_{2}\right) \right) \right) \cdot
_{k}x_{2}\right) \cdot _{k}y_{2} \\[3mm] 
=\left( \gamma \cdot _{h}\rho _{k}\left( r\left( x_{2}\right) \right)
\right) \cdot _{k}\left( x_{2}y_{2}\right) =\left( \gamma \cdot _{h}\rho
_{k}\left( r\left( x_{2}y_{2}\right) \right) \right) \cdot _{k}\left(
x_{2}y_{2}\right) \\[3mm] 
=\gamma \cdot _{kh}\left( x_{2}y_{2}\right)
\end{array}
\end{equation*}

Let $P:G_{1}\rightarrow \mathbf{R}$ be a continuous function with
conditionally compact support such that 
\begin{equation*}
\int P\left( x_{1}\right) d\nu ^{t}\left( x_{1}\right) =1\text{, for all }%
t\in G_{1}^{\left( 0\right) }
\end{equation*}
If $f:\Gamma \star _{kh}G_{2}\rightarrow \mathbf{R}$ is a Borel nonnegative
function, then

$
\begin{array}{l}
\int \int f\left( \gamma ^{-1}\cdot _{kh}x_{2},\gamma ^{-1}\right) d\lambda
^{\rho _{kh}\left( r\left( x_{2}\right) \right) }\left( \gamma \right) d\eta
_{s}\left( x_{2}\right) \\[3mm] 
=\int \int \int f\left( \gamma ^{-1}\cdot _{kh}x_{2},\gamma ^{-1}\right)
d\lambda ^{\rho _{kh}\left( r\left( x_{2}\right) \right) }\left( \gamma
\right) P\left( x_{1}\right) d\nu ^{\rho _{k}\left( r\left( x_{2}\right)
\right) }\left( x_{1}\right) d\eta _{s}\left( x_{2}\right) \text{.}
\end{array}
$

The following sequence of changes of variables 
\begin{equation*}
\begin{array}{l}
\text{1. }\left( x_{2},x_{1}\right) \rightarrow \left( x_{1}^{-1}\cdot
_{k}x_{2},x_{1}\right) \text{(using the quasi-invariance of }\eta _{s}\text{)%
} \text{2. }x_{1}\rightarrow x_{1}^{-1} \\ 
\text{3. }\left( x_{1},\gamma \right) \rightarrow \left( \gamma ^{-1}\cdot
_{h}x_{1},x_{1}\right) \text{(using the quasi-invariance of }\nu _{\rho
_{k}\left( r\left( x_{2}\right) \right) }\text{)} \\ 
\text{4. }x_{1}\rightarrow x_{1}^{-1} \\ 
\text{5. \ }\left( x_{2},x_{1}\right) \rightarrow \left( x_{1}^{-1}\cdot
_{k}x_{2},x_{1}\right) \text{(using the quasi-invariance of }\eta _{s}\text{)%
}
\end{array}
\end{equation*}
transforms the preceding integral into 
\begin{equation*}
\int \int \int f\left( x_{2},\gamma \right) g\left( x_{2},x_{1,}\gamma
\right) P\left( \gamma ^{-1}\cdot _{h}x_{1}\right) d\nu ^{\rho _{k}\left(
r\left( x_{2}\right) \right) }\left( x_{1}\right) d\lambda ^{\rho
_{kh}\left( r\left( x_{2}\right) \right) }\left( \gamma \right) d\eta
_{s}\left( x_{2}\right)
\end{equation*}
where 
\begin{eqnarray*}
g\left( x_{2},x_{1,}\gamma \right) &=&\Delta _{k}\left( x_{1}^{-1}\cdot
_{k}x_{2},\left( \gamma ^{-1}\cdot _{h}x_{1}\right) ^{-1}\right) ^{-1}\Delta
_{h}\left( x_{1},\gamma \right) ^{-1}\Delta _{k}\left( x_{2},x_{1}\right)
^{-1} \\
&=&\left( \Delta _{k}\left( x_{2},x_{1}\right) \Delta _{k}\left(
x_{1}^{-1}\cdot _{k}x_{2},\left( \gamma ^{-1}\cdot _{h}x_{1}\right)
^{-1}\right) \right) ^{-1}\Delta _{h}\left( x_{1},\gamma \right) ^{-1} \\
&=&\Delta _{k}\left( x_{2},\left( \gamma ^{-1}\cdot _{h}r\left( x_{1}\right)
\right) ^{-1}\right) ^{-1}\Delta _{h}\left( x_{1},\gamma \right) ^{-1}\text{.%
}
\end{eqnarray*}

Let us note that for all $\left( x_{1},\gamma \right) \in G_{1}\rtimes
_{h}\Gamma $, and all $x_{2}\in G_{2}$ with $\rho _{k}\left( r\left(
x_{2}\right) \right) =r\left( x_{1}\right) $, $\Delta _{k}\left(
x_{2},\left( \gamma ^{-1}\cdot _{h}r\left( x_{1}\right) \right) \right)
\Delta _{h}\left( x_{1},\gamma \right) $ does not depend on $x_{1}$ but only
on $r\left( x_{1}\right) =\rho _{k}\left( r\left( x_{2}\right) \right) $,
and also it does not depend on $x_{2}$ but only on $r\left( x_{2}\right) $.
For each $\left( x_{2},\gamma \right) \in G_{2}\rtimes _{k}\Gamma $, let us
denote 
\begin{equation*}
\Delta _{kh}\left( x_{2},\gamma \right) =\Delta _{k}\left( x_{2},\left(
\gamma ^{-1}\cdot _{h}r\left( x_{1}\right) \right) ^{-1}\right) \Delta
_{h}\left( x_{1},\gamma \right) \text{.}
\end{equation*}
Consequently, 
\begin{equation*}
\begin{array}{l}
\int \int f\left( \gamma ^{-1}\cdot _{kh}x_{2},\gamma ^{-1}\right) d\lambda
^{\rho _{kh}\left( r\left( x_{2}\right) \right) }\left( \gamma \right) d\eta
_{s}\left( x_{2}\right) \\[3mm] 
=\int \int \int f\left( x_{2},\gamma \right) \Delta _{kh}\left( x_{2},\gamma
\right) ^{-1}P\left( \gamma ^{-1}\cdot _{h}x_{1}\right) d\nu ^{\rho
_{k}\left( r\left( x_{2}\right) \right) }\left( x_{1}\right) d\lambda ^{\rho
_{kh}\left( r\left( x_{2}\right) \right) }\left( \gamma \right) d\eta
_{s}\left( x_{2}\right) \\[3mm] 
=\int \int \int f\left( x_{2},\gamma \right) \Delta _{kh}\left( x_{2},\gamma
\right) ^{-1}P\left( x_{1}\right) d\nu ^{r\left( \gamma ^{-1}\cdot _{h}\rho
_{k}\left( r\left( x_{2}\right) \right) \right) }\left( x_{1}\right)
d\lambda ^{\rho _{kh}\left( r\left( x_{2}\right) \right) }\left( \gamma
\right) d\eta _{s}\left( x_{2}\right) \\[3mm] 
=\int \int \int f\left( x_{2},\gamma \right) \Delta _{kh}\left( x_{2},\gamma
\right) ^{-1}d\lambda ^{\rho _{kh}\left( r\left( x_{2}\right) \right)
}\left( \gamma \right) d\eta _{s}\left( x_{2}\right)
\end{array}
\end{equation*}
Therefore the condition \ref{quasi} in def. \ref{morph} is satisfied if we
take 
\begin{equation*}
\Delta _{kh}\left( x_{2},\gamma \right) =\Delta _{k}\left( x_{2},\left(
\gamma ^{-1}\cdot _{h}r\left( x_{1}\right) \right) \right) \Delta _{h}\left(
x_{1},\gamma \right) ,\left( x_{2},\gamma \right) \in G_{2}\rtimes
_{kh}\Gamma .
\end{equation*}
\end{proof}

\begin{remark}
If $h:\left( \Gamma ,\lambda \right) \mbox{$\,$\rule[0.5ex]{1.1em}{0.2pt}$%
\triangleright\,$} \left( G_{1},\nu \right) $ and $k:\left( G_{1},\lambda
\right) \mbox{$\,$\rule[0.5ex]{1.1em}{0.2pt}$\triangleright\,$} \left(
G_{2},\eta \right) $ are morphisms, then it is easy to see that 
\begin{equation*}
\Delta _{kh}\left( x_{2},\gamma \right) =\Delta _{k}\left( \gamma ^{-1}\cdot
_{kh}x_{2},\gamma ^{-1}\cdot _{h}x_{1}\right) ^{-1}\Delta _{h}\left(
x_{1},\gamma \right) \Delta _{k}\left( x_{2},x_{1}\right)
\end{equation*}
for any $\left( x_{2},\gamma \right) \in G_{2}\rtimes _{kh}\Gamma $ and any $%
x_{1}\in G_{1}^{\rho _{k}\left( r\left( x_{2}\right) \right) }$.
\end{remark}

\begin{proposition}
The class of $\sigma $-compact, lcH-groupoids with the morphisms in the
sense of def. \ref{morph} form a category.
\end{proposition}

\begin{proof}
A straightforward computation shows that the composition of morphisms (in
the sense of def. \ref{comp}) is associative. For each groupoid $\Gamma $
let $l_{\Gamma }$ be the morphism defined in Example \ref{id}. If $h:\Gamma %
\mbox{$\,$\rule[0.5ex]{1.1em}{0.2pt}$\triangleright\,$} G$ and $k:G%
\mbox{$\,$\rule[0.5ex]{1.1em}{0.2pt}$\triangleright\,$} \Gamma $ are
morphisms in the sense of Definition \ref{morph}, then $hl_{\Gamma }=h$ and $%
l_{\Gamma }k=k$.
\end{proof}


\subsubsection{\label{ex_morph}Examples of morphisms}

In this subsection \ we study what becomes a morphism $h:\left( \Gamma
,\lambda \right) \mbox{$\,$\rule[0.5ex]{1.1em}{0.2pt}$\triangleright\,$}
\left( G,\nu \right) $ \ for a particular groupoid $G$. We shall consider
the following cases:

\begin{enumerate}
\item  \textit{Groups}. A group $G$ is a groupoid with $G^{\left( 2\right)
}=G\times G$ and $G^{\left( 0\right) }=\left\{ e\right\} $ (the unit
element).

\item  \textit{Sets}. A set $X$ is a groupoid letting 
\begin{equation*}
X^{\left( 2\right) }=diag\left( X\right) =\left\{ \left( x,x\right) ,x\in
G\right\}
\end{equation*}
and defining the operations by $xx=x$, and $x^{-1}=x$.

Sets and groups are particular cases of group bundles (this means groupoids
for which $r\left( x\right) =d\left( x\right) $ for all $x$).

\item  \textit{Equivalence relations}. Let $\mathcal{E}\subset X\times X$ be
(the graph of) an equivalence relation on the set $X$. Let $\mathcal{E}%
^{\left( 2\right) }=\left\{ \left( \left( x_{1},y_{1}\right) ,\left(
x_{2},y_{2}\right) \right) \in \mathcal{E}\times \mathcal{E}%
:y_{1}=x_{2}\right\} $. With product $\left( x,y\right) \left( y,z\right)
=\left( x,z\right) $ and $\left( x,y\right) ^{-1}=\left( y,x\right) $, $%
\mathcal{E}$ is a principal groupoid. $\mathcal{E}^{\left( 0\right) }$ may
be identified with $X$. Two extreme cases deserve to be single out. If $%
\mathcal{E}=X\times X$, then $\mathcal{E}$ is called the trivial groupoid on 
$X$, while if $\mathcal{E}=diag\left( X\right) $, then $\mathcal{E}$ is
called the co-trivial groupoid on $X$ (and may be identified with the
groupoid in example $2$).

If $G$ is any groupoid, then 
\begin{equation*}
R=\left\{ \left( r\left( x\right) ,d\left( x\right) \right) ,\;x\in G\right\}
\end{equation*}
is an equivalence relation on $G^{\left( 0\right) }$. The groupoid defined
by this equivalence relation is called the \textit{principal groupoid
associated with }$G$.

Any locally compact principal groupoid can be viewed as an equivalence
relation on a locally compact space $X$ having its graph $\mathcal{E}\subset
X\times X$ endowed with a locally compact topology compatible with the
groupoid structure. This topology can be finer than the product topology
induced from $X\times X$. We shall endow the principal groupoid associated
with a groupoid $G$ with the quotient topology induced from $G$ by the map 
\begin{equation*}
\left( r,d\right) :G\rightarrow R,\,\left( r,d\right) \left( x\right)
=\left( r\left( x\right) ,d\left( x\right) \right)
\end{equation*}
This topology consists of the sets whose inverse images by $\left(
r,d\right) $ in $G$ are open.
\end{enumerate}

Let $\Gamma $ and $G$ be two $\sigma $-compact, \emph{lcH}-groupoids,
endowed with the Haar systems $\lambda =\left\{ \lambda ^{u},u\in \Gamma
^{\left( 0\right) }\right\} $, respectively, $\nu =\left\{ \nu ^{t},t\in
G^{\left( 0\right) }\right\} $. Let$\left( \left\{ \beta _{s}^{t}\right\}
_{\left( t,s\right) \in R},\left\{ \tilde{\beta}^{t}\right\} _{\dot{t}\in
G^{\left( 0\right) }/G}\right) $ be the\emph{\ }decomposition of the Haar
system $\left\{ \nu ^{t},t\in G^{\left( 0\right) }\right\} $ over the
principal groupoid associated to $G$ and let $\delta _{G}$ be its associated 
$1$-cocycle.

Let $h:\left( \Gamma ,\lambda \right) \mbox{$\,$\rule[0.5ex]{1.1em}{0.2pt}$%
\triangleright\,$} \left( G,\nu \right) $ be a morphism in the sense of Def. 
\ref{morph}. Let us show that if $G$ is a group bundle , then $\Gamma
|_{\rho _{h}\left( G^{\left( 0\right) }\right) }$ is also a group bundle and
the condition \ref{quasi} in the Def. \ref{morph} is automatically
satisfied. Indeed, let $\gamma \in \Gamma |_{\rho _{h}\left( G^{\left(
0\right) }\right) }$. Then there is $t\in G^{\left( 0\right) }$ such that $%
d\left( \gamma \right) =\rho _{h}\left( t\right) $. We 
\begin{equation*}
r\left( \gamma \right) =\rho _{h}\left( r\left( \gamma \cdot _{h}t\right)
\right) =\rho _{h}\left( d\left( \gamma \cdot _{h}t\right) \right) =\rho
_{h}\left( t\right) =d\left( \gamma \right) \text{.}
\end{equation*}
Therefore $\Gamma |_{\rho _{h}\left( G^{\left( 0\right) }\right) }$ is a
group bundle. Let us prove that the condition \ref{quasi} in the def. \ref
{morph} is automatically satisfied if $G$ is a group bundle. If the
restriction of $r$ to $G^{\prime }$is open, then $\left\{ \beta
_{t}^{t},\,t\in G^{\left( 0\right) }\right\} $ is a Haar system for $%
G^{\prime }$(Lemma 1.3/p. 6 \cite{re2}). In our case $G$ is a group bundle,
consequently, $G=G^{\prime }$. If $\nu =\left\{ \nu ^{t},t\in G^{\left(
0\right) }\right\} $ is a Haar system on $G$, then for each $t\in G^{\left(
0\right) }$, $\nu ^{t}$ is a (left) Haar measure on the locally compact
group $G_{t}^{t}$. By the uniqueness of the Haar measure on $G_{t}^{t}$, it
follows that there is $P\left( t\right) \in \mathbf{R}_{+}^{\ast }$, such
that $\nu ^{t}=P\left( t\right) \beta _{t}^{t}$. Thus the restriction of $%
\delta _{G}$ to $G_{t}^{t}$\ is the modular function for $\nu ^{t}$.
Reasoning in the same way, for each $u\in \rho _{h}\left( G^{\left( 0\right)
}\right) $, $\lambda ^{u}$ is a (left) Haar measure on the locally compact
group $\Gamma _{u}^{u}$, and the restriction of $\delta _{\Gamma }$ to $%
\Gamma _{u}^{u}$\ is the modular function for $\lambda ^{u}$. For each $f\in
C_{c}\left( G\rtimes _{h}\Gamma \right) $ and $t\in G$, we have 
\begin{equation*}
\begin{array}{l}
\int \int f\left( \gamma ^{-1}\cdot _{h}x,\gamma ^{-1}\right) d\lambda
^{\rho _{h}\left( r\left( x\right) \right) }\left( \gamma \right) d\nu
_{t}\left( x\right) \\[3mm] 
=\int \int f\left( \gamma ^{-1}\cdot _{h}x,\gamma ^{-1}\right) d\lambda
^{\rho _{h}\left( t\right) }\left( \gamma \right) d\nu _{t}\left( x\right) \\%
[3mm] 
=\int \int f\left( \gamma ^{-1}\cdot _{h}x,\gamma ^{-1}\right) \delta
_{G}\left( x\right) ^{-1}d\nu ^{t}\left( x\right) d\lambda ^{\rho _{h}\left(
t\right) }\left( \gamma \right) \\[3mm] 
=\int \int f\left( x,\gamma ^{-1}\right) \delta _{G}\left( \left( \gamma
^{-1}\cdot _{h}t\right) ^{-1}x\right) ^{-1}d\nu ^{t}\left( x\right) d\lambda
^{\rho _{h}\left( t\right) }\left( \gamma \right) \\[3mm] 
=\int \int f\left( x,\gamma ^{-1}\right) \delta _{G}\left( \gamma ^{-1}\cdot
_{h}t\right) \delta _{G}\left( x\right) ^{-1}d\nu ^{t}\left( x\right)
d\lambda ^{\rho _{h}\left( t\right) }\left( \gamma \right) \\[3mm] 
=\int \int f\left( x,\gamma \right) \delta _{G}\left( \gamma \cdot
_{h}t\right) \delta _{\Gamma }\left( \gamma ^{-1}\right) \delta _{G}\left(
x\right) ^{-1}d\nu ^{t}\left( x\right) d\lambda ^{\rho _{h}\left( t\right)
}\left( \gamma \right) \\[3mm] 
=\int \int f\left( x,\gamma \right) \delta _{G}\left( \gamma \cdot
_{h}t\right) \delta _{\Gamma }\left( \gamma ^{-1}\right) d\lambda ^{\rho
_{h}\left( t\right) }\left( \gamma \right) d\nu _{t}\left( x\right) \text{.}
\end{array}
\end{equation*}

Hence taking $\Delta _{h}\left( x,\gamma \right) =\delta _{\Gamma }\left(
\gamma \right) \delta _{G}\left( \gamma \cdot _{h}r\left( x\right) \right)
^{-1}$the condition \ref{quasi} in the Definition \ref{morph} is satisfied.
Let us note that if $G$ is a group bundle, then each morphism $h:\left(
\Gamma ,\lambda \right) \mbox{$\,$\rule[0.5ex]{1.1em}{0.2pt}$\triangleright%
\,$} \left( G,\nu \right) $, for which $\rho _{h}:G^{\left( 0\right)
}\rightarrow \Gamma ^{\left( 0\right) }$ is a homeomorphism, can be viewed
as a continuous homomorphism $\varphi $ from $\Gamma $ to $G$ for which the
restriction $\varphi |_{\Gamma ^{\left( 0\right) }}=$ $\varphi ^{\left(
0\right) }:\Gamma ^{\left( 0\right) }\rightarrow G^{\left( 0\right) }$ is a
homeomorphism. Indeed if $\varphi :\Gamma \rightarrow G$ is a groupoid
homomorphism (this means that if $\left( \gamma _{1},\gamma _{2}\right) \in
\Gamma ^{\left( 2\right) }$, then $\left( \varphi \left( \gamma _{1}\right)
,\varphi \left( \gamma _{2}\right) \right) \in G^{\left( 2\right) }$ and $%
\varphi \left( \gamma _{1}\gamma _{2}\right) =\varphi \left( \gamma
_{1}\right) \varphi \left( \gamma _{2}\right) $) and if $\varphi ^{\left(
0\right) }:\Gamma ^{\left( 0\right) }\rightarrow G^{\left( 0\right) }$ is a
homeomorphism, then taking $\rho _{h}=\left( \varphi ^{\left( 0\right)
}\right) ^{-1}:G^{\left( 0\right) }\rightarrow \Gamma ^{\left( 0\right) }$,
and defining 
\begin{equation*}
\gamma \cdot _{h}x=\varphi \left( \gamma \right) x
\end{equation*}
we obtain a morphism in the sense of Definition \ref{morph}. Conversely, if $%
h:\left( \Gamma ,\lambda \right) \mbox{$\,$\rule[0.5ex]{1.1em}{0.2pt}$%
\triangleright\,$} \left( G,\nu \right) $ a morphism in the sense of
Definition \ref{morph}, for which $\rho _{h}:G^{\left( 0\right) }\rightarrow
\Gamma ^{\left( 0\right) }$ is a homeomorphism, then let us define $\varphi
:\Gamma \rightarrow G$ by 
\begin{equation*}
\varphi \left( \gamma \right) =\gamma \cdot _{h}\rho _{h}^{-1}\left( d\left(
\gamma \right) \right) ,\,\gamma \in \Gamma \text{.}
\end{equation*}
If $\left( \gamma _{1},\gamma _{2}\right) \in \Gamma ^{\left( 2\right) }$,
then 
\begin{eqnarray*}
d\left( \varphi \left( \gamma _{1}\right) \right) &=&d\left( \gamma
_{1}\cdot _{h}\rho _{h}^{-1}\left( d\left( \gamma _{1}\right) \right)
\right) =\rho _{h}^{-1}\left( d\left( \gamma _{1}\right) \right) =\rho
_{h}^{-1}\left( r\left( \gamma _{2}\right) \right) =\rho _{h}^{-1}\left(
d\left( \gamma _{2}\right) \right) \\[2mm]
&=&d\left( \gamma _{2}\cdot _{h}\rho _{h}^{-1}\left( d\left( \gamma
_{2}\right) \right) \right) =d\left( \varphi \left( \gamma _{2}\right)
\right) =r\left( \varphi \left( \gamma _{2}\right) \right) .
\end{eqnarray*}
Consequently, $\left( \varphi \left( \gamma _{1}\right) ,\varphi \left(
\gamma _{2}\right) \right) \in G^{\left( 2\right) }$ and 
\begin{eqnarray*}
\varphi \left( \gamma _{1}\gamma _{2}\right) &=&\left( \gamma _{1}\gamma
_{2}\right) \cdot _{h}\rho _{h}^{-1}\left( d\left( \gamma _{2}\right) \right)
\\[2mm]
&=&\gamma _{1}\cdot _{h}\left( \gamma _{2}\cdot _{h}\rho _{h}^{-1}\left(
d\left( \gamma _{2}\right) \right) \right) \\[2mm]
&=&\gamma _{1}\cdot _{h}r\left( \gamma _{2}\cdot _{h}\rho _{h}^{-1}\left(
d\left( \gamma _{2}\right) \right) \right) \left( \gamma _{2}\cdot _{h}\rho
_{h}^{-1}\left( d\left( \gamma _{2}\right) \right) \right) \\[2mm]
&=&\gamma _{1}\cdot _{h}d\left( \gamma _{2}\cdot _{h}\rho _{h}^{-1}\left(
d\left( \gamma _{2}\right) \right) \right) \left( \gamma _{2}\cdot _{h}\rho
_{h}^{-1}\left( d\left( \gamma _{2}\right) \right) \right) \\[2mm]
&=&\left( \gamma _{1}\cdot _{h}\rho _{h}^{-1}\left( d\left( \gamma
_{2}\right) \right) \right) \left( \gamma _{2}\cdot _{h}\rho _{h}^{-1}\left(
d\left( \gamma _{2}\right) \right) \right) \\[2mm]
&=&\left( \gamma _{1}\cdot _{h}\rho _{h}^{-1}\left( d\left( \gamma
_{1}\right) \right) \right) \left( \gamma _{2}\cdot _{h}\rho _{h}^{-1}\left(
d\left( \gamma _{2}\right) \right) \right) \\[2mm]
&=&\varphi \left( \gamma _{1}\right) \varphi \left( \gamma _{2}\right) \text{%
.}
\end{eqnarray*}

The restriction of $\varphi $ to $\Gamma ^{\left( 0\right) }$ is a
homeomorphism, because it coincides with $\rho _{h}^{-1}$.

\qquad Therefore if $\Gamma $ and $G$ are locally compact groups, then the
notion of morphism (cf. Definition \ref{morph}) reduces to the usual notion
of group homomorphism.

If $G$ is a set (see example 2 at beginning of the subsection) and if $%
h:\left( \Gamma ,\lambda \right) \mbox{$\,$\rule[0.5ex]{1.1em}{0.2pt}$%
\triangleright\,$}\left( G,\nu \right) $ is a morphism, then 
\begin{equation*}
\gamma \cdot _{h}x=d\left( \gamma \cdot _{h}x\right) =d\left( x\right) =x
\end{equation*}
for each $\left( \gamma ,x\right) $ with $d\left( \gamma \right) =\rho
_{h}\left( r\left( x\right) \right) =\rho _{h}\left( x\right) $. In this
case a morphism is uniquely determined by the map $\rho _{h}:G\rightarrow
\Gamma ^{\left( 0\right) }$.

Let us now assume that $G\subset $ $X\times X$ \ is an equivalence relation,
where $X$ is locally compact, $\sigma $-compact, Hausdorff space. Let us
endow $G$ with the relative product topology from $X\times X$. Let us also
assume that there is a Haar system on $G$, $\nu =\left\{ \nu ^{t},t\in
G^{\left( 0\right) }\right\} $, and let $\Gamma $ be another groupoid
endowed with the Haar systems $\lambda =\left\{ \lambda ^{u},u\in \Gamma
^{\left( 0\right) }\right\} $. Any morphism in the sense of Definition \ref
{morph}, $h:\left( \Gamma ,\lambda \right) \mbox{$\,$%
\rule[0.5ex]{1.1em}{0.2pt}$\triangleright\,$} \left( G,\nu \right) $,
defines a continuous action of $\Gamma $ on $X$, by 
\begin{equation*}
\gamma \cdot x=r\left( \gamma \cdot _{h}\left( x,x\right) \right) .
\end{equation*}
Conversely, let us consider an action of $\Gamma $ on $X$ with the momentum
map $\rho :X\rightarrow \Gamma ^{\left( 0\right) }$, satisfying $\gamma
\cdot x\symbol{126}x$. Then taking $\rho _{h}=\rho $, and 
\begin{equation*}
\gamma \cdot _{h}\left( x,y\right) =\left( \gamma \cdot x,y\right)
\end{equation*}
we obtain an continuous morphism (in the sense of Definition \ \ref{morph_a}%
). The condition \ref{quasi} is not necessarily satisfied. In order to see
that it is enough to consider $\Gamma =G=X\times X$ (the trivial groupoid on 
$X$ endowed with the product topology). Any Haar system on $X\times X$ is of
the form $\left\{ \varepsilon _{x}\times \mu ,x\in X\right\} $, where $\mu $
is a measure of full support on $X$, and $\varepsilon _{x}$ is the unit
point mass at $x$. If $\left\{ \varepsilon _{x}\times \mu _{1},x\in
X\right\} $ is a Haar system on $\Gamma =X\times X$, and $\left\{
\varepsilon _{x}\times \mu _{2},x\in X\right\} $ is a Haar system on $%
G=X\times X$, then the condition \ref{quasi} in the Definition \ref{morph}
is satisfied if and only if $\mu _{1}$ and $\mu _{2}$ are equivalent measure
(have the same null sets) and the Radon Nikodym derivative is a continuous
function.

Let $G=X\times X$ (the trivial groupoid on $X$ endowed with the product
topology), let $\nu =\left\{ \varepsilon _{x}\times \mu ,x\in X\right\} $ be
a Haar system on $G$. Then any continuous morphism $h$ gives rise to a
continuous action of $\Gamma $ on $X$. Conversely, any continuous action of $%
\Gamma $ on $X$ gives rise to a continuous morphism $h$ from $\Gamma $ to $G$%
. In the hypothesis of Corollary \ref{t}, we can choose the measure $\mu $
such that $h:\left( \Gamma ,\lambda \right) \mbox{$\,$%
\rule[0.5ex]{1.1em}{0.2pt}$\triangleright\,$}\left( G,\nu \right) $ becomes
a morphism in the sense of Definition \ref{morph}.

Let us assume that the associated principal groupoid of $\Gamma $ is proper.
This means that it is a closed subset of $\Gamma ^{\left( 0\right) }\times
\Gamma ^{\left( 0\right) }$ endowed with product topology (or equivalently, $%
\Gamma ^{\left( 0\right) }/\Gamma $ is a Hausdorff space) and the map 
\begin{equation*}
\left( r,d\right) :\Gamma \rightarrow R,\,\left( r,d\right) \left( x\right)
=\left( r\left( x\right) ,d\left( x\right) \right)
\end{equation*}
is an open map when $R$ is endowed with the relative product topology coming
from $\Gamma ^{\left( 0\right) }\times \Gamma ^{\left( 0\right) }$. Let $\mu 
$ be a quasi-invariant measure for the Haar system $\lambda =\left\{ \lambda
^{u},u\in \Gamma ^{\left( 0\right) }\right\} $ on $\Gamma $. It can be shown
that there is a quasi-invariant measure $\mu _{0}$ equivalent to $\mu $ such
that the modular function of $\mu _{0}$ is a continuous function $\delta
_{\Gamma }$. Let $S_{0}$ be the support of $\mu _{0}$. Let us take $X=S_{0}$
and let us consider the action $\Gamma $ on $S_{0}$ defined by $\rho
:S_{0}\rightarrow \Gamma ^{\left( 0\right) },\rho \left( u\right) =u$ for
all $u\in S_{0}$, and $\gamma \cdot d\left( \gamma \right) =r\left( \gamma
\right) $ for all $\gamma \in \Gamma |_{S_{0}}$. It is easy to see that $\mu
_{0}$ is a quasi-invariant measure for the Haar system $\left\{ \varepsilon
_{u}\times \lambda ^{u},u\in S_{0}\right\} $ on $S_{0}\rtimes \Gamma $, and
its modular function is $\delta _{\Gamma }$. Thus we can define a morphism $%
h:\left( \Gamma ,\lambda \right) \mbox{$\,$\rule[0.5ex]{1.1em}{0.2pt}$%
\triangleright\,$} \left( S_{0}\times S_{0},\nu \right) $ in the sense of
Definition \ref{morph} (where $\nu =\left\{ \varepsilon _{u}\times \mu
_{0},u\in S_{0}\right\} $) by

\begin{enumerate}
\item  $\rho _{h}:S_{0}\rightarrow \Gamma ^{\left( 0\right) },\rho \left(
u\right) =u$ for all $u\in S_{0}$.

\item  $\gamma \cdot _{h}\left( u,v\right) =\left( r\left( \gamma \right)
,v\right) $
\end{enumerate}


\section{\label{app}Morphisms on a groupoid $\Gamma $ and the convolution
algebra $C_{c}\left( \Gamma \right) $}

Let $(\Gamma ,\lambda )$ and $(G,\nu )$ be two $\sigma $-compact, \emph{lcH}%
-groupoids with Haar systems.

We associate to each morphism $h:\left( \Gamma ,\lambda \right) %
\mbox{$\,$\rule[0.5ex]{1.1em}{0.2pt}$\triangleright\,$}\left( G,\nu \right) $
an application $\hat{h}$ defined on $C_{c}\left( \Gamma \right) $ in the
following way. For any $f\in C_{c}\left( \Gamma \right) $, 
\begin{equation*}
\hat{h}\left( f\right) :C_{c}\left( G\right) \rightarrow C_{c}\left(
G\right) .
\end{equation*}
is defined by 
\begin{equation*}
\hat{h}\left( f\right) \left( \xi \right) \left( x\right) =\int f\left(
\gamma \right) \xi \left( \gamma ^{-1}\cdot _{h}x\right) \Delta _{h}\left(
x,\gamma \right) ^{-1/2}d\lambda ^{\rho _{h}\left( r\left( x\right) \right)
}\left( \gamma \right)
\end{equation*}
Using a standard argument (\cite{co} 2.2, \cite{re1} II.1) we can prove that 
$\hat{h}\left( f\right) \left( \xi \right) \in C_{c}\left( G\right) $ for
any $\xi \in C_{c}\left( G\right) $. That is, since $G\rtimes _{h}\Gamma $
is a closed subset of the normal space $G\times \Gamma $, the function 
\begin{equation*}
\left( x,\gamma \right) \rightarrow \xi \left( \gamma ^{-1}\cdot
_{h}x\right) \Delta _{h}\left( x,\gamma \right) ^{-1/2}\text{,}
\end{equation*}
may be extended to a bounded continuous function $F$ on $G\times \Gamma $. A
compactness argument shows that for each $\varepsilon >0$ and each $x_{0}\in
G$%
\begin{equation*}
\left\{ x\in G:\,\left| F\left( x,\gamma \right) -F\left( x_{0},\gamma
\right) \right| <\varepsilon \text{ for all }\gamma \in supp\left( f\right)
\right\}
\end{equation*}
is an open subset of $G$ which contains $x_{0}$. Therefore the function 
\begin{equation*}
x\rightarrow F_{x}\,\left[ :G\rightarrow C_{c}\left( G\right) \right]
\end{equation*}
where $F_{x}\left( y\right) =f\left( \gamma \right) F\left( x,\gamma \right) 
$, is continuous. Consequently, 
\begin{equation*}
\left( x,u\right) \rightarrow \int f\left( \gamma \right) F\left( x,\gamma
\right) d\lambda ^{u}\left( \gamma \right) \,\,\left[ :G\times \Gamma
^{\left( 0\right) }\rightarrow \mathbf{C}\right]
\end{equation*}
is a continuous function, and so is the function 
\begin{equation*}
x\rightarrow \int f\left( \gamma \right) F\left( x,\gamma \right) d\lambda
^{\rho _{k}\left( r\left( x\right) \right) }\left( \gamma \right) \,\ \, 
\left[ :G\rightarrow \mathbf{C}\right]
\end{equation*}
(being its composition with $x\rightarrow \left( x,\rho _{h}\left( r\left(
x\right) \right) \right) $.

\begin{lemma}
\label{dens} Let $h:\left( \Gamma ,\lambda \right) \mbox{$\,$%
\rule[0.5ex]{1.1em}{0.2pt}$\triangleright\,$} \left( G,\nu \right) $ be a
morphism of $\sigma $-compact, lcH-groupoids with Haar systems and let $\hat{%
h}$ be the application defined above . Then 
\begin{equation*}
\left\{ \hat{h}\left( f\right) \xi :f\in C_{c}\left( \Gamma \right) \text{, }%
\xi \in C_{c}\left( G\right) \right\}
\end{equation*}
is dense in $C_{c}\left( G\right) $ with the inductive limit topology.
\end{lemma}

\begin{proof}
We shall use a similar argument as Jean Renault used in proof of prop.
II.1.9/p. 56 \cite{re1}. Since $\Gamma ^{\left( 0\right) }$ is a paracompact
space, it follows that there is \ a fundamental system of $d$-relatively
compact neighborhood $\left\{ U_{\alpha }\right\} _{\alpha }$ of $\Gamma
^{\left( 0\right) }$. Let $U_{0}$ be a $d$-relatively compact neighborhood
of $\Gamma ^{\left( 0\right) }$ such that $U_{\alpha }\subset U_{0}$ for all 
$\alpha $. Let $\left\{ K_{\alpha }\right\} _{\alpha }$ be a net of compact
subsets of $\Gamma ^{\left( 0\right) }$ increasing to $\Gamma ^{\left(
0\right) }$. \ Let $e_{\alpha }\in C_{c}\left( \Gamma \right) $ be a
nonnegative function such that 
\begin{eqnarray*}
supp\left( e_{\alpha }\right) &\subset &U_{\alpha } \\
\int e_{\alpha }\left( \gamma \right) d\lambda ^{u}\left( \gamma \right) &=&1%
\text{ \ for all }u\in K_{\alpha }\text{.}
\end{eqnarray*}
We claim that for any $\xi \in C_{c}\left( G\right) $, $\left\{ \hat{h}%
\left( e_{\alpha }\right) \xi \right\} _{\alpha }$ converges to $\xi $ in
the inductive limit topology. Let $\xi \in C_{c}\left( G\right) $ and $%
\varepsilon >0$. Let $K$ be the support of $\xi $. Then 
\begin{eqnarray*}
U_{0}\cdot _{h}K &=&\left\{ \gamma \cdot _{h}x:\gamma \in U_{0},\,x\in
K,\,\,r\left( \gamma \right) =\rho _{h}\left( r\left( x\right) \right)
\right\} \\[3mm]
&=&\left( U_{0}\cap d^{-1}\left\{ \rho _{h}\left( r\left( K\right) \right)
\right\} \right) \cdot _{h}K
\end{eqnarray*}
is a compact subset of $G$. A compactness argument shows that 
\begin{equation*}
W_{\varepsilon }=\left\{ \gamma \in \Gamma :\,\left| \xi \left( \gamma
^{-1}\cdot _{h}x\right) -\xi \left( x\right) \right| <\varepsilon \text{ for
all }x\in U_{0}\cdot _{h}K,\rho _{h}\left( r\left( x\right) \right) =r\left(
\gamma \right) \right\}
\end{equation*}
is an open subset of $\Gamma $ which contains $\Gamma ^{\left( 0\right) }$.
If $\gamma \in W_{\varepsilon }\cap U_{0}$, then 
\begin{equation*}
\,\left| \xi \left( \gamma ^{-1}\cdot _{h}x\right) -\xi \left( x\right)
\right| <\varepsilon \text{ for all }x\text{ satisfying }\rho _{h}\left(
r\left( x\right) \right) =r\left( \gamma \right)
\end{equation*}
(because if $x\notin U_{0}\cdot _{h}K$, then $x$ and $\gamma ^{-1}\cdot
_{h}x\notin K=supp\left( \xi \right) $, and hence $\xi \left( \gamma
^{-1}\cdot _{h}x\right) =\xi \left( x\right) =0$). Since $\Delta _{h}$ is a
continuous function and a homomorphism from $G\times _{h}\Gamma $ to $%
\mathbf{R}_{+}^{\ast }$, it follows that there exist an open neighborhood $%
L_{\varepsilon }$ of $\Gamma ^{\left( 0\right) }$ such that 
\begin{equation*}
\left| \Delta _{h}^{-1/2}\left( x,\gamma \right) -1\right| <\varepsilon
\end{equation*}
for all $\left( x,\gamma \right) \in \left( K\times L_{\varepsilon }\right)
\cap G\times _{h}\Gamma $. Then for any $\alpha $ such that $U_{\alpha
}\subset W_{\varepsilon }\cap L_{\varepsilon }$ and $\rho _{h}\left( r\left(
U_{0}\cdot _{h}K\right) \right) \subset K_{\alpha }$, $supp\left( \hat{h}%
\left( e_{\alpha }\right) \xi \right) $ is contained in $U_{0}\cdot _{h}K$.
For all $x\in U_{0}\cdot _{h}K$ we have 
\begin{eqnarray*}
&&\left| \hat{h}\left( e_{\alpha }\right) \xi \left( x\right) -\xi \left(
x\right) \right| \\[3mm]
&=&\left| \int e_{\alpha }\left( \gamma \right) \xi \left( \gamma ^{-1}\cdot
_{h}x\right) \Delta _{h}\left( x,\gamma \right) ^{-1/2}-e_{\alpha }\left(
\gamma \right) \xi \left( x\right) d\lambda ^{\rho _{h}\left( r\left(
x\right) \right) }\left( \gamma \right) \right| \\[3mm]
&\leq &\int e_{\alpha }\left( \gamma \right) \left| \xi \left( \gamma
^{-1}\cdot _{h}x\right) -\xi \left( x\right) \right| \Delta _{h}\left(
x,\gamma \right) ^{-1/2}d\lambda ^{\rho _{h}\left( r\left( x\right) \right)
}\left( \gamma \right) + \\
&&\,\ \ \ \ \ \ \ \ \ \ \ \ \ \ \ \ \ +\left| \xi \left( x\right) \right|
\int e_{\alpha }\left( \gamma \right) \left| \Delta _{h}^{-1/2}\left(
x,\gamma \right) -1\right| d\lambda ^{\rho _{h}\left( r\left( x\right)
\right) }\left( \gamma \right) \\[3mm]
&\leq &2\varepsilon +\sup\limits_{x}\left| \xi \left( x\right) \right|
\varepsilon \text{.}
\end{eqnarray*}
Thus $\left| \hat{h}\left( e_{\alpha }\right) \xi -\xi \right| $ converges
to $0$ in the inductive limit topology.
\end{proof}

\begin{proposition}
\label{ac} Let $h:\left( \Gamma ,\lambda \right) \mbox{$\,$%
\rule[0.5ex]{1.1em}{0.2pt}$\triangleright\,$} \left( G_{1},\nu \right) $ and 
$k:\left( G_{1},\lambda \right) \mbox{$\,$\rule[0.5ex]{1.1em}{0.2pt}$%
\triangleright\,$} \left( G_{2},\eta \right) $ be morphisms of $\sigma $%
-compact lcH-groupoids with Haar systems. Then 
\begin{equation*}
\hat{k}\left( \hat{h}\left( f\right) \xi _{1}\right) \xi _{2}=\hat{kh}\left(
f\right) \left( \hat{k}\left( \xi _{1}\right) \xi _{2}\right)
\end{equation*}
for all $f\in C_{c}\left( \Gamma \right) $, $\xi _{1}\in C_{c}\left(
G_{1}\right) $ and $\xi _{2}\in C_{c}\left( G_{2}\right) $.
\end{proposition}

\begin{proof}
Let $f\in C_{c}\left( \Gamma \right) $, $\xi _{1}\in C_{c}\left(
G_{1}\right) $ and $\xi _{2}\in C_{c}\left( G_{2}\right) $. For all $\left(
x_{2},x_{1},\gamma \right) \in G_{2}\times G_{1}\times \Gamma $, such that $%
\left( x_{1},\gamma \right) \in G_{1}\rtimes _{h}\Gamma $ and $\left(
x_{2},x_{1}\right) \in G_{2}\rtimes _{h}G_{1}$ let us denote 
\begin{eqnarray*}
F\left( x_{2},x_{1},\gamma \right) &=&\xi _{1}\left( \gamma ^{-1}\cdot
_{h}x_{1}\right) \xi _{2}\left( x_{1}^{-1}\cdot _{k}x_{2}\right) \\
g\left( x_{2},x_{1},\gamma \right) &=&\Delta _{h}\left( x_{1},\gamma \right)
^{-1/2}\Delta _{k}\left( x_{2},x_{1}\right) ^{-1/2}.
\end{eqnarray*}
Then we have 
\begin{equation*}
\begin{array}{l}
\hat{k}\left( \hat{h}\left( f\right) \xi _{1}\right) \xi _{2}\left(
x_{2}\right) \\[3mm] 
=\int \int f\left( \gamma \right) F\left( x_{2},x_{1},\gamma \right) g\left(
x_{2},x_{1},\gamma \right) \lambda ^{\rho _{h}\left( r\left( x_{1}\right)
\right) }\left( \gamma \right) d\nu ^{\rho _{k}\left( r\left( x_{2}\right)
\right) }\left( x_{1}\right) \\[3mm] 
=\int \int f\left( \gamma \right) F\left( x_{2},x_{1},\gamma \right) g\left(
x_{2},x_{1},\gamma \right) d\nu ^{\rho \left( r\left( x_{2}\right) \right)
}\left( x_{1}\right) d\lambda ^{\rho _{kh}\left( r\left( x_{2}\right)
\right) }\left( \gamma \right) \\[3mm] 
=\int \int f\left( \gamma \right) F_{1}\left( x_{2},x_{1},\gamma \right)
g_{1}\left( x_{2},x_{1},\gamma \right) d\nu ^{r\left( \gamma ^{-1}\cdot
_{h}\rho _{k}\left( r\left( x_{2}\right) \right) \right) }\left(
x_{1}\right) d\lambda ^{\rho _{kh}\left( r\left( x_{2}\right) \right)
}\left( \gamma \right) \text{,}
\end{array}
\end{equation*}
where 
\begin{eqnarray*}
F_{1}\left( x_{2},x_{1},\gamma \right) &=&F\left( x_{2},\left( \gamma
^{-1}\cdot _{h}\rho _{k}\left( r\left( x_{2}\right) \right) \right)
^{-1}x_{1},\gamma \right) \\
&=&\xi _{1}\left( x_{1}\right) \xi _{2}\left( \left( \left( \gamma
^{-1}\cdot _{h}\rho _{k}\left( r\left( x_{2}\right) \right) \right)
^{-1}x_{1}\right) ^{-1}\cdot _{k}x_{2}\right) \\
&=&\xi _{1}\left( x_{1}\right) \xi _{2}\left( x_{1}^{-1}\left( \gamma
^{-1}\cdot _{h}\rho _{k}\left( r\left( x_{2}\right) \right) \right) \cdot
_{k}x_{2}\right) \\[2mm]
&=&\xi _{1}\left( x_{1}\right) \xi _{2}\left( x_{1}^{-1}\cdot _{k}\left(
\gamma ^{-1}\cdot _{kh}x_{2}\right) \right)
\end{eqnarray*}
and 
\begin{equation*}
\begin{array}{l}
g_{1}\left( x_{2},x_{1},\gamma \right) =g\left( x_{2},\left( \gamma
^{-1}\cdot _{h}\rho _{k}\left( r\left( x_{2}\right) \right) \right)
^{-1}x_{1},\gamma \right) \\[3mm] 
=\Delta _{h}\left( \left( \gamma ^{-1}\cdot _{h}\rho _{k}\left( r\left(
x_{2}\right) \right) \right) ^{-1}x_{1},\gamma \right) ^{-1/2}\Delta
_{k}\left( x_{2},\left( \gamma ^{-1}\cdot _{h}\rho _{k}\left( r\left(
x_{2}\right) \right) \right) ^{-1}x_{1}\right) ^{-1/2} \\ 
=\Delta _{h}\left( \rho _{k}\left( r\left( x_{2}\right) \right) ,\gamma
\right) ^{-1/2}\Delta _{k}\left( x_{2},\left( \gamma ^{-1}\cdot _{h}\rho
_{k}\left( r\left( x_{2}\right) \right) \right) ^{-1}x_{1}\right) ^{-1/2} \\%
[3mm] 
=\Delta _{h}\left( \rho _{k}\left( r\left( x_{2}\right) \right) ,\gamma
\right) ^{-1/2}\Delta _{k}\left( x_{2},\left( \gamma ^{-1}\cdot _{h}\rho
_{k}\left( r\left( x_{2}\right) \right) \right) ^{-1}\right) ^{-1/2}\Delta
_{k}\left( \gamma ^{-1}\cdot _{kh}x_{2},x_{1}\right) ^{-1/2} \\[3mm] 
=\Delta _{kh}\left( x_{2},\gamma \right) ^{-1/2}\Delta _{k}\left(
x_{2},\left( \gamma ^{-1}\cdot _{h}\rho _{k}\left( r\left( x_{2}\right)
\right) \right) ^{-1}\right) ^{-1/2}\Delta _{k}\left( \gamma ^{-1}\cdot
_{kh}x_{2},x_{1}\right) ^{-1/2}
\end{array}
\end{equation*}

Consequently, 
\begin{equation*}
\hat{k}\left( \hat{h}\left( f\right) \xi _{1}\right) \xi _{2}\left(
x_{2}\right) =\hat{kh}\left( f\right) \left( \hat{k}\left( \xi _{1}\right)
\xi _{2}\right) \left( x_{2}\right) \text{, for all }x_{2}\in G_{2\text{.}}
\end{equation*}
\end{proof}

For any locally compact, second countable, Hausdorff groupoid $G$ endowed
with a Haar system $\nu =\left\{ \nu ^{t},t\in G^{\left( 0\right) }\right\} $%
, $C_{c}\left( G\right) $ is an algebra under convolution of function. For $%
f $, $g\in C_{c}\left( G\right) $ the convolution is defined by:

\begin{equation*}
f\ast g\left( x\right) =\int f\left( y\right) g\left( y^{-1}x\right) d\nu
^{r\left( x\right) }\left( y\right)
\end{equation*}
and the involution by 
\begin{equation*}
f^{\ast }\left( x\right) =\overline{f\left( x^{-1}\right) }\text{.}
\end{equation*}
Moreover, under these operations, $C_{c}\left( G\right) $ becomes a
topological $\ast $-algebra. Let us note that the involutive algebraic
structure on $C_{c}\left( G\right) $ defined above depends on the Haar
system $\nu =\left\{ \nu ^{t},t\in G^{\left( 0\right) }\right\} $. When it
will be necessary to emphasis the role of $\nu $ in this structure, we shall
write $C_{c}\left( G,\nu \right) $.

It is easy to see that for any $f$, $g\in C_{c}\left( G,\nu \right) $%
\begin{equation*}
f\ast g=\hat{l}\left( f\right) g\text{, }
\end{equation*}
where $l:\left( G,\nu \right) \mbox{$\,$\rule[0.5ex]{1.1em}{0.2pt}$%
\triangleright\,$} \left( G,\nu \right) $ is the morphism defined in Example 
\ref{id}\ : $\rho _{l}=id_{G^{\left( 0\right) }}$ and $x\cdot _{l}y=xy$
(multiplication on $G$)$.$

For each $f\in C_{c}\left( G\right) $, let us denote by $\left\| f\right\|
_{I}$ the maximum of $\sup\limits_{t}\int \left| f\left( x\right) \right|
d\nu ^{t}\left( x\right) $ and $\sup\limits_{t}\int \left| f\left( x\right)
\right| d\nu _{t}\left( x\right) $. A straightforward computation shows that 
$\left\| \cdot \right\| _{I}$ is a norm on $C_{c}\left( G\right) $ and 
\begin{eqnarray*}
\left\| f\right\| _{I} &=&\left\| f^{\ast }\right\| _{I} \\
\left\| f\ast g\right\| _{I} &\leq &\left\| f\right\| _{I}\left\| g\right\|
_{I}
\end{eqnarray*}
for all $f,\,g\in C_{c}\left( G\right) $.

\begin{proposition}
\label{herm} Let $h:\left( \Gamma ,\lambda \right) \mbox{$\,$%
\rule[0.5ex]{1.1em}{0.2pt}$\triangleright\,$} \left( G,\nu \right) $ be a
morphism of $\sigma$-compact lcH-groupoids with Haar systems. Then 
\begin{equation*}
\xi _{2}^{\ast }\ast \left( \hat{h}\left( f\right) \xi _{1}\right) =\left( 
\hat{h}\left( f^{\ast }\right) \xi _{2}\right) ^{\ast }\ast \xi _{1}
\end{equation*}
for all $f\in C_{c}\left( \Gamma \right) $ and $\xi _{1},\,\xi _{2}\in
C_{c}\left( G,\nu \right) $.
\end{proposition}

\begin{proof}
If $f\in C_{c}\left( \Gamma \right) $ and $\xi _{1},\,\xi _{2}\in
C_{c}\left( G\right) $, then 
\begin{equation*}
\begin{array}{l}
\xi _{2}^{\ast }\ast \left( \hat{h}\left( f\right) \xi _{1}\right) \left(
x\right) = \\[3mm] 
=\int \xi _{2}^{\ast }\left( y\right) \int f\left( \gamma \right) \xi
_{1}\left( \gamma ^{-1}\cdot _{h}\left( y^{-1}x\right) \right) \Delta
_{h}\left( y^{-1}x,\gamma \right) ^{-1/2}d\lambda ^{\rho _{h}\left( r\left(
y^{-1}x\right) \right) }\left( \gamma \right) d\nu ^{r\left( x\right)
}\left( y\right) \\[3mm] 
=\int \overline{\xi _{2}\left( y^{-1}\right) }\int f\left( \gamma \right)
\xi _{1}\left( \gamma ^{-1}\cdot _{h}\left( y^{-1}x\right) \right) \Delta
_{h}\left( y^{-1},\gamma \right) ^{-1/2}d\lambda ^{\rho _{h}\left( r\left(
y^{-1}x\right) \right) }\left( \gamma \right) d\nu ^{r\left( x\right)
}\left( y\right) \\[3mm] 
=\int \int \overline{\xi _{2}\left( y\right) }f\left( \gamma \right) \xi
_{1}\left( \gamma ^{-1}\cdot _{h}\left( yx\right) \right) \Delta _{h}\left(
y,\gamma \right) ^{-1/2}d\lambda ^{\rho _{h}\left( r\left( yx\right) \right)
}\left( \gamma \right) d\nu _{r\left( x\right) }\left( y\right)
\end{array}
\end{equation*}

The change of variable $\left( y,\,\gamma \right) \rightarrow \left( \gamma
^{-1}\cdot _{h}y,\,\gamma ^{-1}\right) $ transforms the preceding integral
into 
\begin{equation*}
\begin{array}{l}
=\int \int \overline{\xi _{2}\left( \gamma ^{-1}\cdot _{h}y\right) }f\left(
\gamma ^{-1}\right) \Delta _{h}\left( y,\gamma \right) ^{-1/2}d\lambda
^{\rho _{h}\left( r\left( y\right) \right) }\left( \gamma \right) \xi
_{1}\left( yx\right) d\nu _{r\left( x\right) }\left( y\right) \\[3mm] 
=\int \int \overline{\xi _{2}\left( \gamma ^{-1}\cdot _{h}y^{-1}\right) }%
f\left( \gamma ^{-1}\right) \Delta _{h}\left( y^{-1},\gamma \right)
^{-1/2}d\lambda ^{\rho _{h}\left( d\left( y\right) \right) }\left( \gamma
\right) \xi _{1}\left( y^{-1}x\right) d\nu ^{r\left( x\right) }\left(
y\right) \\[3mm] 
=\int \overline{\int f^{\ast }\left( \gamma \right) \xi _{2}\left( \gamma
^{-1}\cdot _{h}y^{-1}\right) \Delta _{h}\left( y^{-1},\,\gamma \right)
^{-1/2}d\lambda ^{\rho _{h}\left( d\left( y\right) \right) }\left( \gamma
\right) }\xi _{1}\left( y^{-1}x\right) d\nu ^{r\left( x\right) }\left(
y\right) \\[3mm] 
=\left( \hat{h}\left( f^{\ast }\right) \xi _{2}\right) ^{\ast }\ast \xi _{1}
\end{array}
\end{equation*}
\end{proof}


\section{Representations associated to morphisms}

\begin{proposition}
\label{rep} Let $h:\left( \Gamma ,\lambda \right) \mbox{$\,$%
\rule[0.5ex]{1.1em}{0.2pt}$\triangleright\,$} \left( G,\nu \right) $ be a
morphism of $\sigma $-compact, lcH-groupoids with Haar systems. For $t\in
G^{\left( 0\right) }$ and $f\in C_{c}\left( \Gamma \right) $, let us define
the operator $\pi _{h,t}\left( f\right) :\;L^{2}\left( G,\,\nu _{t}\right) %
\mbox{$\,$\rule[0.5ex]{1.1em}{0.2pt}$\triangleright\,$} L^{2}\left( G,\,\nu
_{t}\right) $ by 
\begin{equation*}
\pi _{h,t}\left( f\right) \xi \left( x\right) =\int f\left( \gamma \right)
\xi \left( \gamma ^{-1}\cdot _{h}x\right) \Delta _{h}\left( x,\gamma \right)
^{-1/2}d\lambda ^{\rho _{h}\left( r\left( x\right) \right) }\left( \gamma
\right)
\end{equation*}
for all $\xi \in L^{2}\left( G,\,\nu _{t}\right) $ and $x\in G$.Then for any 
$f\in C_{c}\left( \Gamma \right) $ 
\begin{equation*}
\left\| \pi _{h,t}\left( f\right) \right\| \leq \left\| f\right\| _{I}\text{,%
}
\end{equation*}
and $\pi _{h,t}$ is a representation of $C_{c}\left( \Gamma ,\,\lambda
\right) $ (a $\ast $-homomorphism from $C_{c}\left( \Gamma ,\,\lambda
\right) $ into $\emph{B}\left( L^{2}\left( G,\,\nu _{t}\right) \right) $,
that is continuous with respect to the inductive limit topology on $%
C_{c}\left( \Gamma \right) $ and the weak operator topology on $\emph{B}%
\left( L^{2}\left( G,\,\nu _{t}\right) \right) $).
\end{proposition}

\begin{proof}
If $f\in C_{c}\left( \Gamma \right) $, $\xi ,\zeta \in L^{2}\left( G,\,\nu
_{t}\right) $, then 
\begin{equation*}
\begin{array}{l}
\left| \left\langle \pi _{h,t}\left( f\right) \xi ,\zeta \right\rangle
\right| = \\[3mm] 
=\left| \int \int f\left( \gamma \right) \xi \left( \gamma ^{-1}\cdot
_{h}x\right) \Delta _{h}\left( x,\gamma \right) ^{-1/2}d\lambda ^{\rho
_{h}\left( r\left( x\right) \right) }\left( \gamma \right) \overline{\zeta
\left( x\right) }d\nu _{t}\left( x\right) \right| \\[3mm] 
\leq \int \int \left| f\left( \gamma \right) \right| \left| \xi \left(
\gamma ^{-1}\cdot _{h}x\right) \right| \left| \zeta \left( x\right) \right|
\Delta _{h}\left( x,\gamma \right) ^{-1/2}d\lambda ^{\rho _{h}\left( r\left(
x\right) \right) }\left( \gamma \right) d\nu _{t}\left( x\right) \\[3mm] 
\leq \left( \int \int \left| f\left( \gamma \right) \right| \left| \xi
\left( \gamma ^{-1}\cdot _{h}x\right) \right| ^{2}\Delta _{h}\left( x,\gamma
\right) ^{-1}d\lambda ^{\rho _{h}\left( r\left( x\right) \right) }\left(
\gamma \right) d\nu _{t}\left( x\right) \right) ^{1/2}\cdot \\ 
\,\ \ \ \ \ \ \ \ \ \ \ \ \ \ \ \ \ \ \ \ \ \ \ \ \ \ \ \ \ \ \cdot \left(
\int \int \left| f\left( \gamma \right) \right| \left| \zeta \left( x\right)
\right| ^{2}d\lambda ^{\rho _{h}\left( r\left( x\right) \right) }\left(
\gamma \right) d\nu _{t}\left( x\right) \right) ^{1/2} \\[3mm] 
=\left( \int \int \left| f\left( \gamma ^{-1}\right) \right| \left| \xi
\left( x\right) \right| ^{2}d\lambda ^{\rho _{h}\left( r\left( x\right)
\right) }\left( \gamma \right) d\nu _{t}\left( x\right) \right) ^{1/2}\cdot
\\ 
\ \,\ \ \ \ \ \ \ \ \ \ \ \ \ \ \ \ \ \ \ \ \ \ \ \ \ \ \ \ \cdot \left(
\int \int \left| f\left( \gamma \right) \right| d\lambda ^{\rho _{h}\left(
r\left( x\right) \right) }\left( \gamma \right) \left| \zeta \left( x\right)
\right| ^{2}d\nu _{t}\left( x\right) \right) ^{1/2} \\[3mm] 
\leq \left\| f\right\| _{I}\left\| \xi \right\| _{2}\left\| \zeta \right\|
_{2}\text{.}
\end{array}
\end{equation*}
Thus $\left\| \pi _{h,t}\left( f\right) \right\| \leq \left\| f\right\| _{I}$
for any $f\in C_{c}\left( \Gamma \right) $. Let us prove that $\pi
_{h,t}:C_{c}\left( \Gamma \right) \mbox{$\,$\rule[0.5ex]{1.1em}{0.2pt}$%
\triangleright\,$} \emph{B}\left( L^{2}\left( G,\,\nu _{t}\right) \right) $
is a $\ast $-homomorphism. Let $f\in C_{c}\left( \Gamma \right) $, $\xi
,\zeta \in L^{2}\left( G,\,\nu _{t}\right) $. We have 
\begin{equation*}
\begin{array}{l}
\left\langle \pi _{h,t}\left( f\right) \xi ,\zeta \right\rangle = \\[3mm] 
=\int \int f\left( \gamma \right) \xi \left( \gamma ^{-1}\cdot _{h}x\right)
\Delta _{h}\left( x,\gamma \right) ^{-1/2}d\lambda ^{\rho _{h}\left( r\left(
x\right) \right) }\left( \gamma \right) \overline{\zeta \left( x\right) }%
d\nu _{t}\left( x\right) \\[3mm] 
=\int \int f\left( \gamma ^{-1}\right) \xi \left( x\right) \overline{\zeta
\left( \gamma ^{-1}\cdot _{h}x\right) }\Delta _{h}\left( x,\gamma \right)
^{-1/2}d\lambda ^{\rho _{h}\left( r\left( x\right) \right) }\left( \gamma
\right) d\nu _{t}\left( x\right) \\[3mm] 
=\overline{\int \int f^{\ast }\left( \gamma \right) \zeta \left( \gamma
^{-1}\cdot _{h}x\right) \Delta _{h}\left( x,\gamma \right) ^{-1/2}\overline{%
\xi \left( x\right) }d\lambda ^{\rho _{h}\left( r\left( x\right) \right)
}\left( \gamma \right) d\nu _{t}\left( x\right) } \\[3mm] 
=\overline{\left\langle \pi _{h,t}\left( f^{\ast }\right) \zeta ,\xi
\right\rangle }\text{.}
\end{array}
\end{equation*}
Hence $\pi _{h,t}\left( f\right) ^{\ast }=\pi _{h,t}\left( f^{\ast }\right) $
for all $f\in C_{c}\left( \Gamma \right) $. If $f,g\in C_{c}\left( \Gamma
\right) $ and $\xi \in L^{2}\left( G,\,\nu _{t}\right) $, we have 
\begin{equation*}
\begin{array}{l}
\pi _{h,t}\left( f\ast g\right) \xi \left( x\right) = \\[3mm] 
=\int f\ast g\left( \gamma \right) \xi \left( \gamma ^{-1}\cdot _{h}x\right)
\Delta _{h}\left( x,\gamma \right) ^{-1/2}d\lambda ^{\rho _{h}\left( r\left(
x\right) \right) }\left( \gamma \right) \\[3mm] 
=\int \int f\left( \gamma _{1}\right) g\left( \gamma _{1}^{-1}\gamma \right)
d\lambda ^{r\left( \gamma \right) }\left( \gamma _{1}\right) \xi \left(
\gamma ^{-1}\cdot _{h}x\right) \Delta _{h}\left( x,\gamma \right)
^{-1/2}d\lambda ^{\rho _{h}\left( r\left( x\right) \right) }\left( \gamma
\right) \\[3mm] 
=\int \int f\left( \gamma _{1}\right) g\left( \gamma _{1}^{-1}\gamma \right)
\xi \left( \gamma ^{-1}\cdot _{h}x\right) \Delta _{h}\left( x,\gamma \right)
^{-1/2}d\lambda ^{\rho _{h}\left( r\left( x\right) \right) }\left( \gamma
_{1}\right) d\lambda ^{\rho _{h}\left( r\left( x\right) \right) }\left(
\gamma \right) \\[3mm] 
=\int \int f\left( \gamma _{1}\right) g\left( \gamma \right) \xi \left(
\left( \gamma _{1}\gamma \right) ^{-1}\cdot _{h}x\right) \Delta _{h}\left(
x,\gamma _{1}\gamma \right) ^{-1/2}d\lambda ^{d\left( \gamma _{1}\right)
}\left( \gamma \right) d\lambda ^{\rho _{h}\left( r\left( x\right) \right)
}\left( \gamma _{1}\right) \\[3mm] 
=\int f\left( \gamma _{1}\right) \int g\left( \gamma \right) \xi \left(
\gamma ^{-1}\cdot _{h}\left( \gamma _{1}^{-1}\cdot _{h}x\right) \right)
\Delta _{h}\left( \gamma _{1}^{-1}\cdot _{h}x,\gamma \right) ^{-1/2}d\lambda
^{d\left( \gamma _{1}\right) }\left( \gamma \right) \cdot \\ 
\,\ \ \ \ \ \ \ \ \ \ \ \ \ \ \ \ \ \ \ \ \ \ \ \ \ \ \ \ \ \ \ \ \ \ \ \ \
\ \ \ \ \ \ \ \ \ \ \ \ \ \ \ \ \ \ \ \ \ \ \ \ \cdot \Delta _{h}\left(
x,\gamma _{1}\right) ^{-1/2}d\lambda ^{\rho _{h}\left( r\left( x\right)
\right) }\left( \gamma _{1}\right) \\[3mm] 
=\pi _{h,t}\left( f\right) \left( \pi _{h,t}\left( g\right) \xi \right)
\left( x\right) \text{,}
\end{array}
\end{equation*}
for all $x\in G$. Thus $\pi _{h,t}\left( f\ast g\right) =\pi _{h,t}\left(
f\right) \pi _{h,t}\left( g\right) $ for all $f,g\in C_{c}\left( \Gamma
\right) $.
\end{proof}

\begin{remark}
The representation $\pi _{h,t}:C_{c}\left( \Gamma \right) \rightarrow \emph{B%
}\left( L^{2}\left( G,\,\nu _{t}\right) \right) $ defined in the preceding
proposition is non-degenerate in the sense that 
\begin{equation*}
\left\{ \pi _{h,t}\left( f\right) \xi :\,f\in C_{c}\left( \Gamma \right)
,\xi \in L^{2}\left( G,\,\nu _{t}\right) \right\}
\end{equation*}
is dense in $L^{2}\left( G,\,\nu _{t}\right) $. For any $f\in C_{c}\left(
\Gamma \right) $ and $\xi \in C_{c}\left( G\right) $, $\pi _{h,t}\left(
f\right) \xi =\hat{h}\left( f\right) \xi $, where $\hat{h}$ is the
application introduced in Section \ref{app}. Indeed, by Lemma \ref{dens} 
\begin{equation*}
\left\{ \pi _{h,t}\left( f\right) \xi :\,f\in C_{c}\left( \Gamma \right)
,\xi \in C_{c}\left( G\right) \right\}
\end{equation*}
is dense in $C_{c}\left( G\right) $ with the inductive limit topology and a
fortiori with the $L^{2}\left( G,\,\nu _{t}\right) $ topology.
\end{remark}

\begin{proposition}
\label{rep_comp} Let $h:\left( \Gamma ,\lambda \right) \mbox{$\,$%
\rule[0.5ex]{1.1em}{0.2pt}$\triangleright\,$} \left( G_{1},\nu \right) $ and 
$k:\left( G_{1},\lambda \right) \mbox{$\,$\rule[0.5ex]{1.1em}{0.2pt}$%
\triangleright\,$} \left( G_{2},\eta \right) $ be morphisms of $\sigma$%
-compact lcH groupoids with Haar systems. Let $\hat{h}$ the application
associated to $h$ introduced in the Section \ref{app}, and $\pi _{k,s}$
(respectively, $\pi _{kh,s}$) be the representation associated to $k$
(respectively, $kh$) defined in Proposition \ref{rep}. Then 
\begin{equation*}
\pi _{k,s}\left( \hat{h}\left( f\right) \xi _{1}\right) \xi _{2}=\pi
_{kh,s}\left( f\right) \pi _{k,s}\left( \xi _{1}\right) \xi _{2}
\end{equation*}
for all $s\in G_{2}^{\left( 0\right) }$, $f\in C_{c}\left( \Gamma \right) $, 
$\xi _{1}\in C_{c}\left( G_{1}\right) $ and $\xi _{2}\in L^{2}\left(
G_{2},\eta _{s}\right) $.
\end{proposition}

\begin{proof}
It follows using the same computation as in the proof of Proposition \ref{ac}%
.
\end{proof}


\section{A $C^{\ast }$-algebra associated to a locally compact $\protect%
\sigma $-compact groupoid}

In this section we define a $C^{\ast }$-algebra associated to a locally
compact $\sigma $-compact groupoid. The construction is similar to that made
in Section 5 \cite{st}. Like in \cite{st} we shall show that $\left( \Gamma
,\lambda \right) \mbox{$\,$\rule[0.5ex]{1.1em}{0.2pt}$\triangleright\,$}
C^{\ast }\left( \Gamma ,\lambda \right) $ is a covariant functor from the
category of locally compact, $\sigma $-compact, Hausdorff groupoids endowed
with Haar systems to the category of $C^{\ast }$-algebras. \ \ 

\begin{definition}
\label{norm}Let $(\Gamma ,\lambda )$ be a $\sigma $-compact, lcH-groupoid
endowed with a Haar and let $f\in C_{c}\left( \Gamma \right) $. For any
morphism $h:\left( \Gamma ,\lambda \right) \mbox{$\,$%
\rule[0.5ex]{1.1em}{0.2pt}$\triangleright\,$}\left( G,\nu \right) $
\thinspace (where $G$ is a $\sigma $-compact, lcH-groupoid endowed with a
Haar system) let us define 
\begin{equation*}
\left\| f\right\| _{h}=\sup\limits_{t}\left\| \pi _{h,t}\left( f\right)
\right\|
\end{equation*}
where $\pi _{h,t}$ is the representation associated to $h$ and $t$ defined
in prop. \ref{rep}. Let us also define 
\begin{equation*}
\left\| f\right\| =\sup\limits_{t}\left\| f\right\| _{h}
\end{equation*}
where $h$ runs over all morphism defined on $\left( \Gamma ,\lambda \right) $%
.
\end{definition}

\begin{remark}
Clearly $\left\| \cdot \right\| =\sup\limits_{h}\left\| \cdot \right\| _{h}$
defined above is a $C^{\ast }$-semi-norm. Let $l:\left( \Gamma ,\lambda
\right) \mbox{$\,$\rule[0.5ex]{1.1em}{0.2pt}$\triangleright\,$}\left( \Gamma
,\lambda \right) $ be the morphism defined in Example \ref{id}:$\rho
_{l}=id_{\Gamma ^{0}}$ and $\gamma \cdot _{l}x=\gamma x$ (multiplication on $%
\Gamma $). Then the representation associated to $l$ and $u\in \Gamma
^{\left( 0\right) }$, $\pi _{h,u}:C_{c}\left( \Gamma \right) \rightarrow 
\emph{B}\left( L^{2}\left( G,\,\,\lambda _{u}\right) \right) $, is defined
by 
\begin{equation*}
\pi _{l,u}\left( f\right) \xi \left( x\right) =\int f\left( \gamma \right)
\xi \left( \gamma ^{-1}x\right) d\lambda ^{r\left( x\right) }\left( \gamma
\right) =f\ast \xi \left( x\right)
\end{equation*}
by for all $f\in C_{c}\left( \Gamma \right) $ and $\xi \in L^{2}\left(
G,\,\,\lambda _{u}\right) $. Therefore for all $f\in C_{c}\left( \Gamma
\right) $ 
\begin{equation*}
\left\| f\right\| _{l}=\left\| f\right\| _{red}\,\text{,}
\end{equation*}
the reduced norm of $f$ (def. 2.36/p. 50 \cite{mu}, or def. II.2.8./p. 82 
\cite{re1}). According to prop. II.1.11/p. 58 \cite{re1}, $\left\{ \pi
_{l,u}\right\} _{u}$ is a faithful family of representations of $C_{c}\left(
\Gamma ,\lambda \right) $, so $\left\| \cdot \right\| _{red}$ is a norm.
Hence $\left\| \cdot \right\| =sup_{h}\left\| \cdot \right\| _{h}$ (where $h$
runs over all morphism defined on $\left( \Gamma ,\lambda \right) $) is a
norm on $C_{c}\left( \Gamma ,\lambda \right) $.
\end{remark}

\begin{definition}
\label{alg}Let $(\Gamma,\lambda) $ be a $\sigma $-compact, lcH-groupoid
endowed with Haar system. The $C^{\ast }$-algebra $C^{\ast }\left( \Gamma
,\lambda \right) $ is defined to be the completion of $C_{c}\left( \Gamma
,\lambda \right) $ in the norm $\left\| \cdot \right\| =sup_{h}\left\| \cdot
\right\| _{h}$, where $h$ runs over all morphism defined on $\left( \Gamma
,\lambda \right) $.
\end{definition}

\begin{remark}
Let $\Gamma $ be\ a locally compact, second countable, Hausdorff groupoid
endowed with Haar system $\lambda =\left\{ \lambda ^{u},u\in \Gamma ^{\left(
0\right) }\right\} $. For $f\in C_{c}\left( \Gamma ,\lambda \right) $,\ the
full $C^{\ast }$-norm is defined by 
\begin{equation*}
\left\| f\right\| _{full}=\sup\limits_{L}\left\| L\left( f\right) \right\|
\end{equation*}
where $L$ is a non-degenerate representation of $\ C_{c}\left( \Gamma
,\lambda \right) $, i.e. a $\ast $-homomorphism from $C_{c}\left( \Gamma
,\,\lambda \right) $ into $\emph{B}\left( H\right) $, for some Hilbert space 
$H$, that is continuous with respect to the inductive limit topology on $%
C_{c}\left( \Gamma \right) $ and the weak operator topology on $\emph{B}%
\left( H\right) $, and is such that the linear span of 
\begin{equation*}
\left\{ L\left( g\right) \xi :\,g\in C_{c}\left( \Gamma \right) ,\,\xi \in
H\right\}
\end{equation*}
is dense in $H$. We have the following inequalities 
\begin{equation*}
\left\| f\right\| _{red}\leq \left\| f\right\| \leq \left\| f\right\| _{full}
\end{equation*}
for all $f\in C_{c}\left( \Gamma ,\lambda \right) $, where $\left\| \cdot
\right\| =\sup\limits_{h}\left\| \cdot \right\| _{h}$ is the norm introduced
in Definition \ref{norm}.The full $C^{\ast }$algebra $C_{full}^{\ast }\left(
\Gamma ,\lambda \right) $ and the reduced $C^{\ast }$algebra $C_{red}^{\ast
}\left( \Gamma ,\lambda \right) $ are defined respectively as the completion
of the algebra $C_{c}\left( \Gamma ,\lambda \right) $ for the full norm $%
\left\| \cdot \right\| _{full}$, and the reduced norm $\left\| \cdot
\right\| _{red}$. According to Proposition 6.1.8/p.146 \cite{ar}, if $\left(
\Gamma ,\lambda \right) $ is measurewise amenable (Definition 3.3.1/p. 82 
\cite{ar}), then $C_{full}^{\ast }\left( \Gamma ,\lambda \right)
=C_{red}^{\ast }\left( \Gamma ,\lambda \right) $. Thus if $\left( \Gamma
,\lambda \right) $ is measurewise amenable, then $C^{\ast }\left( \Gamma
,\lambda \right) $\ (the $C^{\ast }$-algebra introduced in Definition \ref
{alg}), $C_{full}^{\ast }\left( \Gamma ,\lambda \right) $ and $C_{red}^{\ast
}\left( \Gamma ,\lambda \right) $ coincide.
\end{remark}

\begin{notation}
Let $\Gamma $ be\ a locally compact, second countable, Hausdorff groupoid
endowed with Haar system $\lambda =\left\{ \lambda ^{u},u\in \Gamma ^{\left(
0\right) }\right\} $. Let $\mu $ be a quasi-invariant measure . Let $\Delta
\mu $ be the modular function associated to $\left\{ \lambda ^{u},u\in
\Gamma ^{\left( 0\right) }\right\} $ and $\mu $. Let $\lambda _{1}$ be the
measure induced by $\mu $ on $\Gamma $, and $\lambda _{0}=\Delta _{\mu }^{-%
\frac{1}{2}}\lambda _{1}$.

For $f\in L^{1}\left( G,\lambda _{0}\right) $ let us define$\qquad $%
\begin{equation*}
\left\| f\right\| _{II,\mu }=\sup \left\{ \int |f\left( \gamma \right)
j\left( d\left( \gamma \right) k\left( r\left( \gamma \right) \right)
\right) |d\lambda _{0}\left( \gamma \right) \right\}
\end{equation*}
the supremum being taken over all $\ j,k\in L^{2}\left( \Gamma ^{\left(
0\right) },\mu \right) $with $\int |j|^{2}d\mu =\int |k|^{2}d\mu =1.$

Let 
\begin{equation*}
II\left( G,\lambda ,\mu \right) =\left\{ f\in L^{1}\left( G,\lambda
_{0}\right) ,\left\| f\right\| _{II,\mu }<\infty \right\} \text{.}
\end{equation*}
$II\left( G,\lambda ,\mu \right) $ is a Banach $\ast $-algebra under the ,
convolution 
\begin{equation*}
f\ast g\left( \gamma _{1}\right) =\int f\left( \gamma \right) g\left( \gamma
^{-1}\gamma _{1}\right) d\lambda ^{r\left( \gamma _{1}\right) }\left( \gamma
\right)
\end{equation*}
and the involution 
\begin{equation*}
f^{\ast }\left( \gamma \right) =\overline{f\left( \gamma ^{-1}\right) }\text{%
.}
\end{equation*}
$C_{c}\left( \Gamma ,\lambda \right) $ is a $\ast $-subalgebra of $II\left(
G,\lambda ,\mu \right) $ for any quasi-invariant measure $\mu $. If $\mu
_{1} $ and $\mu _{2}$ are two equivalent quasi-invariant measures, then $%
\left\| f\right\| _{II,\mu _{1}}=\left\| f\right\| _{II,\mu _{2}}$ for all $%
f $. Let us denote 
\begin{equation*}
\left\| f\right\| _{II}=\sup \left\{ \left\| f\right\| _{II,\mu }\right\}
\end{equation*}
the supremum being taken over all quasi-invariant measure $\mu $ on $\Gamma
^{\left( 0\right) }$, and let us note that it is enough to consider one
quasi-invariant measure in each class.

For any $f\in II\left( G,\lambda ,\mu \right) $, 
\begin{equation*}
k\rightarrow \left( u\rightarrow \int f\left( \gamma \right) k\left( d\left(
\gamma \right) \right) d\lambda ^{u}\left( \gamma \right) \right) \text{.}
\end{equation*}
is a bounded operator $II_{\mu }\left( f\right) $ on $L^{2}\left( G^{\left(
0\right) },\mu \right) $ and $\left\| II_{\mu }\left( \left| f\right|
\right) \right\| =\left\| f\right\| _{II,\mu }$. Moreover, $f\rightarrow
II_{\mu }\left( f\right) $ is a norm-decreasing $\ast $-representation $%
II_{\mu }$\ of $II\left( G,\lambda ,\mu \right) $. The restriction of $%
II_{\mu }$ to $C_{c}\left( \Gamma ,\lambda \right) $ is a representation of $%
C_{c}\left( \Gamma ,\,\lambda \right) $ called the trivial representation on 
$\mu $.

Every representation $\left( \mu ,\Gamma ^{\left( 0\right) }\ast \mathcal{H}%
,L\right) $(see Definition 3.20/p.68 \cite{mu}) of $\Gamma $ can be
integrated into a representation, still denoted by $L$, of $II\left(
G,\lambda ,\mu \right) $. The relation between the two representation is: 
\begin{equation*}
\left\langle L\left( f\right) \xi _{1},\xi _{2}\right\rangle =\int f\left(
\gamma \right) \left\langle L\left( \gamma \right) \xi _{1}\left( d\left(
\gamma \right) \right) ,\xi _{2}\left( r\left( \gamma \right) \right)
\right\rangle \lambda _{0}\left( \gamma \right)
\end{equation*}
where $f\in $ $II\left( G,\lambda ,\mu \right) $, $\xi _{1},\xi _{2}\in
\int_{G^{\left( 0\right) }}^{\oplus }\mathcal{H}\left( u\right) d\mu \left(
u\right) $. Conversely, every non-degenerate $\ast $-representation of any
suitably large $\ast $-algebra of $II\left( G,\lambda ,\mu \right) $ (in
particular, $C_{c}\left( \Gamma ,\lambda \right) $) is equivalent to a
representation obtained this fashion (see Section 3 \cite{ha2}, Proposition
II.1.17/p. 52\cite{re1}, Proposition 4.2 \cite{re3} or Proposition 3.23/p.
70, Theorem 3.29/p. 74 \cite{mu}). If $L$ is the integrated form of a
representation, $\left( \mu ,G^{\left( 0\right) }\ast \mathcal{H},L\right) $%
, of the groupoid $G$, then 
\begin{equation*}
\left| \left\langle L\left( f\right) \xi ,\eta \right\rangle \right| \leq
\left\langle II_{\mu }\left( \left| f\right| \right) \widetilde{\xi },%
\widetilde{\eta }\right\rangle
\end{equation*}
where $\widetilde{\xi }\left( u\right) =\left\| \xi \left( u\right) \right\|
.\,$Therefore $\left\| L\left( f\right) \right\| \leq \left\| II\mu \left(
\left| f\right| \right) \right\| =\left\| f\right\| _{II,\mu }$ $\leq
\left\| f\right\| _{II}$.

If $f\in C_{c}\left( \Gamma ,\lambda \right) $, then 
\begin{equation*}
\left\| \left| f\right| \right\| _{II,\mu }=\left\| II_{\mu }\left( \left|
f\right| \right) \right\| \leq \left\| \left( \left| f\right| \right)
\right\| _{full}\leq \left\| \left| f\right| \right\| _{II}.
\end{equation*}
Thus if $f\in C_{c}\left( \Gamma ,\lambda \right) $ and $f\geq 0$, then $%
\left\| f\right\| _{II}=\left\| f\right\| _{full}$.
\end{notation}

\begin{proposition}
Let $\Gamma $ be\ a locally compact, second countable, Hausdorff groupoid
endowed with Haar system $\lambda =\left\{ \lambda ^{u},u\in \Gamma ^{\left(
0\right) }\right\} $. Let $\left\| \cdot \right\| $ be the norm on the $%
C^{\ast }$-algebra introduced in Definition \ref{alg}, and $\left\| \cdot
\right\| _{full}$ the norm on the full $C^{\ast }$algebra associated with $%
\Gamma $ and $\lambda $ introduced by Renault in \cite{re1}. Let us assume
the principal associated groupoid of $\Gamma $ is a proper groupoid. Then
for any quasi invariant measure $\mu $ on $\Gamma ^{\left( 0\right) }$ and
any $f\in C_{c}\left( \Gamma ,\lambda \right) $ 
\begin{equation*}
\left\| II_{\mu }\left( f\right) \right\| \leq \left\| f\right\| \text{,}
\end{equation*}
and if $f\geq 0$, then $\left\| f\right\| =\left\| f\right\| _{full}=\left\|
II\left( f\right) \right\| $.
\end{proposition}

\begin{proof}
For each $f\in C_{c}\left( \Gamma ,\lambda \right) $ we have $\left\|
f\right\| \leq \left\| f\right\| _{full}$. If $f\in C_{c}\left( \Gamma
,\lambda \right) $ and $f\geq 0$, then $\left\| f\right\| _{full}=\left\|
f\right\| _{II}=\sup \left\{ \left\| f\right\| _{II,\mu }\right\} $, the
supremum being taken over all quasi-invariant measure $\mu $ on $\Gamma
^{\left( 0\right) }$. Let $\mu $ be a quasi-invariant measure. We have shown
at the end of Subsection \ref{ex_morph}\ that if the principal associated
groupoid of $\Gamma $ is a proper groupoid, then there is a quasi-invariant
measure $\mu _{0}$ equivalent to $\mu $ such that the modular function of $%
\mu _{0}$ is a continuous function $\delta _{\Gamma }$. Let $S_{0}$ be the
support of $\mu _{0}$. Let us consider the action $\Gamma $ on $S_{0}$
defined by $\rho :S_{0}\rightarrow \Gamma ^{\left( 0\right) },$ $\rho \left(
u\right) =u$ for all $u\in S_{0}$, and $\gamma \cdot d\left( \gamma \right)
=r\left( \gamma \right) $ for all $\gamma \in \Gamma |_{S_{0}}$. It is easy
to see that $\mu _{0}$ is a quasi-invariant measure for the Haar system $\nu
=\left\{ \varepsilon _{u}\times \lambda ^{u},u\in S_{0}\right\} $ on $%
S_{0}\rtimes \Gamma $, and its modular function is $\delta _{\Gamma }$. Thus
we can define a morphism $h_{\mu }:\left( \Gamma ,\lambda \right)
\rightarrow \left( S_{0}\times S_{0},\nu \right) $ in the sense of
Definition \ref{morph} (where $\nu =\left\{ \varepsilon _{u}\times \mu
_{0},u\in S_{0}\right\} $) by

\begin{enumerate}
\item  $\rho _{h}:S_{0}\rightarrow \Gamma ^{\left( 0\right) },\rho \left(
u\right) =u$ for all $u\in S_{0}$.

\item  $\gamma \cdot _{h_{\mu }}\left( u,v\right) =\left( r\left( \gamma
\right) ,v\right) $
\end{enumerate}

The representations $\pi _{h_{\mu },u}$ associated to the morphism $h_{\mu
}:\left( \Gamma ,\lambda \right) \mbox{$\,$\rule[0.5ex]{1.1em}{0.2pt}$%
\triangleright\,$}\left( S_{0}\times S_{0},\nu \right) $ as in Proposition 
\ref{rep}, can be identified with the trivial representation $II_{\mu }$ of $%
C_{c}\left( \Gamma ,\lambda \right) $. Hence for any $f\in C_{c}\left(
\Gamma ,\lambda \right) $%
\begin{equation*}
\left\| f\right\| \geq \left\| f\right\| _{h_{\mu }}=\left\| II_{\mu }\left(
f\right) \right\| \text{.}
\end{equation*}
Therefore $f\in C_{c}\left( \Gamma ,\lambda \right) $ and $f\geq 0$, then 
\begin{equation*}
\left\| f\right\| _{full}=\left\| f\right\| _{II}=\sup_{\mu }\left\|
f\right\| _{II,\mu }=\sup_{\mu }\left\| II_{\mu }\left( f\right) \right\|
=\sup_{\mu }\left\| f\right\| _{h_{\mu }}\geq \left\| f\right\| .
\end{equation*}
\end{proof}

The following propositions are slightly modified version of props.5.2/p. 27
and 5.3/p. 27 \cite{st}.

\begin{proposition}
Let $h:\left( \Gamma ,\lambda \right) \mbox{$\,$\rule[0.5ex]{1.1em}{0.2pt}$%
\triangleright\,$} \left( G,\nu \right) $ be a morphism of $\sigma $%
-compact, lcH-groupoids with Haar systems and let $\hat{h}$ be the mapping
defined in Section \ref{app}.Then $\hat{h}$ extends to $\ast $-homomorphism 
\begin{equation*}
C^{\ast }\left( h\right) :C^{\ast }\left( \Gamma ,\lambda \right)
\rightarrow M\left( C^{\ast }\left( G,\nu \right) \right) \text{,}
\end{equation*}
where $M\left( C^{\ast }\left( G,\nu \right) \right) $ is multiplier algebra
of $C^{\ast }\left( G,\nu \right) $, with the property that $C^{\ast }\left(
h\right) \left( C^{\ast }\left( \Gamma ,\lambda \right) \right) C^{\ast
}\left( G,\nu \right) $ is dense in $C^{\ast }\left( G,\nu \right) $.
\end{proposition}

\begin{proof}
Let $G_{1}$ be a locally compact, $\sigma $-compact, Hausdorff groupoid
endowed with a Haar system $\eta =\left\{ \eta ^{s},\,s\in G_{1}^{\left(
0\right) }\right\} $ and let $k:\left( G,\nu \right) \mbox{$\,$%
\rule[0.5ex]{1.1em}{0.2pt}$\triangleright\,$} \left( G_{1},\eta \right) $ be
a morphism. Let $\left\{ \pi _{k,s}\right\} _{s}$ be the family of
representations defined in Proposition \ref{rep}. $\ $According to
Proposition \ref{rep_comp}, for all $s\in G_{1}^{\left( 0\right) }$, $f\in
C_{c}\left( \Gamma \right) $, $\xi _{1}\in C_{c}\left( G\right) $ and $\xi
_{2}\in L^{2}\left( G_{1},\eta _{s}\right) $, 
\begin{equation*}
\pi _{k,s}\left( \hat{h}\left( f\right) \xi _{1}\right) \xi _{2}=\pi
_{kh,s}\left( f\right) \pi _{k,s}\left( \xi _{1}\right) \xi _{2}\text{.}
\end{equation*}
Hence for all $s\in G_{1}^{\left( 0\right) }$, $f\in C_{c}\left( \Gamma
\right) $ and $\xi \in C_{c}\left( G\right) $, 
\begin{equation*}
\left\| \pi _{k,s}\left( \hat{h}\left( f\right) \xi \right) \right\|
=\left\| \pi _{kh,s}\left( f\right) \right\| \left\| \pi _{k,s}\left( \xi
\right) \right\| \leq \left\| f\right\| _{C^{\ast }\left( \Gamma ,\lambda
\right) }\left\| \xi \right\| _{C^{\ast }\left( G,\nu \right) }
\end{equation*}

Let $\left\| \cdot \right\| _{C^{\ast }\left( \Gamma ,\lambda \right) }$
(and respectively, $\left\| \cdot \right\| _{C^{\ast }\left( G,\nu \right) }$%
) be the norm on the $C^{\ast }$-algebra introduced in Definition \ref{alg}.
Thus 
\begin{equation*}
\left\| \hat{h}\left( f\right) \xi \right\| _{C^{\ast }\left( G,\nu \right)
}\leq \left\| f\right\| _{C^{\ast }\left( \Gamma ,\lambda \right) }\left\|
\xi \right\| _{C^{\ast }\left( G,\nu \right) }
\end{equation*}
Since $C_{c}\left( G\right) $ is dense in $C^{\ast }\left( G,\nu \right) $,
it follows that \ $\hat{h}\left( f\right) $ extends to bounded linear map $%
C^{\ast }\left( h\right) \left( f\right) $ on $C^{\ast }\left( G,\nu \right) 
$. By Proposition \ref{herm}, for all $f\in C_{c}\left( \Gamma \right) $ and 
$\xi _{1},\,\xi _{2}\in C_{c}\left( G,\nu \right) $, 
\begin{equation*}
\xi _{2}^{\ast }\ast \left( \hat{h}\left( f\right) \xi _{1}\right) =\left( 
\hat{h}\left( f^{\ast }\right) \xi _{2}\right) ^{\ast }\ast \xi _{1}\text{.}
\end{equation*}
Using the density of $C_{c}\left( G\right) $ in $C^{\ast }\left( G,\nu
\right) $ and the continuity of $C^{\ast }\left( h\right) \left( f\right) $
we have 
\begin{equation*}
\xi _{2}^{\ast }\ast \left( C^{\ast }\left( h\right) \left( f\right) \xi
_{1}\right) =\left( C^{\ast }\left( h\right) \left( f^{\ast }\right) \xi
_{2}\right) ^{\ast }\ast \xi _{1}\text{,}
\end{equation*}
for all $\xi _{1},\,\xi _{2}\in C^{\ast }\left( G,\nu \right) $. Hence $%
C^{\ast }\left( h\right) \left( f\right) $ admits a Hermitian adjoint $%
C^{\ast }\left( h\right) \left( f^{\ast }\right) $,\thinspace and therefore $%
C^{\ast }\left( h\right) \left( f\right) \in M\left( C^{\ast }\left( G,\nu
\right) \right) $. Since $C_{c}\left( \Gamma \right) $ is dense in $C^{\ast
}\left( \Gamma ,\lambda \right) $, it follows that $C^{\ast }\left( h\right) 
$\ extends to $C^{\ast }\left( \Gamma ,\lambda \right) $.

By Lemma \ref{dens} 
\begin{equation*}
\left\{ \hat{h}\left( f\right) \xi :f\in C_{c}\left( \Gamma \right) \text{, }%
\xi \in C_{c}\left( G\right) \right\}
\end{equation*}
is dense in $C_{c}\left( G\right) $ with the inductive limit topology and a
fortiori in $C^{\ast }\left( G,\nu \right) $. Therefore $C^{\ast }\left(
h\right) \left( C^{\ast }\left( \Gamma ,\lambda \right) \right) C^{\ast
}\left( G,\nu \right) $ is dense in $C^{\ast }\left( G,\nu \right) .$
\end{proof}

\begin{proposition}
Let $\Gamma $, $G_{1}$ and $G_{2}$ be locally compact, $\sigma $-compact,
Hausdorff groupoids. Let $\lambda =\left\{ \lambda ^{u},u\in \Gamma ^{\left(
0\right) }\right\} $ (respectively, $\nu =\left\{ \nu ^{t},t\in
G_{1}^{\left( 0\right) }\right\} $, $\eta =\left\{ \eta ^{s},s\in
G_{2}^{\left( 0\right) }\right\} $) be a Haar system on $\Gamma $
(respectively, on $G_{1}$, $G_{2}$). Let $h:\left( \Gamma ,\lambda \right) %
\mbox{$\,$\rule[0.5ex]{1.1em}{0.2pt}$\triangleright\,$} \left( G_{1},\nu
\right) $ and $k:\left( G_{1},\lambda \right) \mbox{$\,$%
\rule[0.5ex]{1.1em}{0.2pt}$\triangleright\,$} \left( G_{2},\eta \right) $ be
morphisms. Then 
\begin{equation*}
C^{\ast }\left( kh\right) =C^{\ast }\left( k\right) C^{\ast }\left( h\right) 
\text{.}
\end{equation*}
\end{proposition}

\begin{proof}
Let us denote by $\hat{C}^{\ast }\left( k\right) :$ $M\left( C^{\ast }\left(
G_{1},\nu \right) \right) \mbox{$\,$\rule[0.5ex]{1.1em}{0.2pt}$%
\triangleright\,$} M\left( C^{\ast }\left( G_{2},\eta \right) \right) $ the
unique extension of $C^{\ast }\left( k\right) $. In order to show that $%
C^{\ast }\left( kh\right) =\hat{C}^{\ast }\left( k\right) C^{\ast }\left(
h\right) $, it is enough to prove that $C^{\ast }\left( kh\right) \left(
f\right) =\hat{C}^{\ast }\left( k\right) \left( C^{\ast }\left( h\right)
\right) \left( f\right) $ for $f\in C_{c}\left( \Gamma \right) $. As a
consequence of Lemma \ref{dens} 
\begin{equation*}
\left\{ \hat{k}\left( \xi _{1}\right) \xi _{2}:\xi _{1}\in C_{c}\left(
G_{1}\right) ,\,\xi _{2}\in C_{c}\left( G_{2}\right) \right\}
\end{equation*}
is dense in $C^{\ast }\left( G_{2},\eta \right) $. Thus for proving $C^{\ast
}\left( kh\right) \left( f\right) =\hat{C}^{\ast }\left( k\right) C^{\ast
}\left( h\right) \left( f\right) $ it is enough to prove 
\begin{eqnarray*}
\hat{kh}\left( f\right) \left( \hat{k}\left( \xi _{1}\right) \xi _{2}\right)
&=&\hat{C}^{\ast }\left( k\right) \left( \hat{h}\left( f\right) \right) \hat{%
k}\left( \xi _{1}\right) \xi _{2} \\
&=&\hat{k}\left( \hat{h}\left( f\right) \xi _{1}\right) \xi _{2}
\end{eqnarray*}
but this is true (see Proposition \ref{ac}).
\end{proof}

\begin{remark}
We have constructed a covariant functor $\left( \Gamma ,\lambda \right)
\rightarrow C^{\ast }\left( \Gamma ,\lambda \right) $, $h\rightarrow C^{\ast
}\left( h\right) $ from the category of locally compact, $\sigma $-compact,
Hausdorff groupoids endowed with Haar systems to the category of $C^{\ast }$%
-algebras (in the sense of \cite{wo}). The hypothesis of $\sigma $%
-compactness is not really necessary. It is enough to work with groupoids
which are locally compact and normal, and for which the unit spaces have
conditionally-compact neighborhoods (for instance, paracompact unit spaces).
\end{remark}

\textbf{Acknowledgements.}M. Buneci was partly supported by the MEdC-ANCS grant ET65/2005 and by the Postdoctoral Training Program HPRN-CT-2002-0277. P. Stachura was supported by by Polish KBN grant 115/E-343/SPB/6.PR UE/DIE50/2005-2008.

\bigskip

{\small M\u{a}d\u{a}lina Buneci}

{\small Department of Automatics}

{\small University Constantin Br\^{a}ncu\c{s}i of T\^{a}rgu-Jiu}

{\small Bul. Republicii 1, 210152, T\^{a}rgu-Jiu, Romania}

{\small e-mail: ada@utgjiu.ro}

\bigskip

{\small Piotr Stachura}

{\small Department of Mathematical Methods in Physics}

University of {\small Warsaw }

ul. {\small Ho\.{z}a\ 74,\ 00-682,\ Warszawa,\ Poland}

{\small e-mail: stachura@fuw.edu.pl}

\end{document}